\documentclass[onefignum,onetabnum,sort&compress]{siamonline190516}

\usepackage{lipsum}
\usepackage{graphicx}
\usepackage{epstopdf}
\usepackage{algorithmic}
\ifpdf
  \DeclareGraphicsExtensions{.eps,.pdf,.png,.jpg}
\else
  \DeclareGraphicsExtensions{.eps}
\fi

\headers{On the Koopman operator of algorithms}{Felix Dietrich, Thomas N. Thiem, and Ioannis G. Kevrekidis}

\title{On the Koopman operator of algorithms\thanks{Submitted to the editors on \today.\funding{This work was partially funded by the DARPA/Lagrange program (Drs. Fahroo and Lewis) and by the ARO through a MURI (Drs. S. Stanton and M. Munson).}}}

\author{Felix Dietrich\thanks{Department of Chemical and Biomolecular Engineering and Department of Applied Mathematics and Statistics, Johns Hopkins University and Department of Informatics, Technical University of Munich
  (\email{felix.dietrich@tum.de}).}
\and Thomas N. Thiem\thanks{Department of Chemical and Biomolecular Engineering, Princeton University
  (\email{tthiem@princeton.edu}).}
\and Ioannis
Kevrekidis\thanks{Department of Chemical and Biomolecular Engineering and Department of Applied Mathematics and Statistics, Johns Hopkins University and JHMI 
  (\email{yannis@princeton.edu}).}}

\ifpdf
\hypersetup{
  pdftitle={On the Koopman operator of algorithms},
  pdfauthor={Felix Dietrich, Thomas N. Thiem, Ioannis G. Kevrekidis}
}
\fi

\usepackage[utf8]{inputenc}
\usepackage{amsmath}
\usepackage{amsfonts}
\usepackage{amssymb}
\usepackage{tikz}

\usepackage{soul}

\usepackage{color}

\usepackage{graphicx}
\graphicspath{{./figures/}{./}}

\newcommand{\figureThirdWidth}{.3\textwidth}
\newcommand{\figureHalfWidth}{4.5cm}
\newcommand{\figureFullWidth}{.9\textwidth}

\newcommand{\K}     {\mathcal{K}} 
\newcommand{\A}     {A_{\K}}      
\newcommand{\Keval} {\lambda}     
\newcommand{\Aeval} {\omega}      
\newcommand{\Kefunc}{\phi}        
\newcommand{\flow}  {{S}}         
\newcommand{\obs}   {g}           
\newcommand{\figref}[1]{figure~\ref{fig:#1}}
\newcommand{\Figref}[1]{Figure~\ref{fig:#1}}
\newcommand{\secref}[1]{section~\ref{sec:#1}}

\newcommand{\eqnref}[1]{equation~\eqref{eq:#1}}
\newcommand{\Eqnref}[1]{Equation~\eqref{eq:#1}}

\newcommand{\Real}{\text{Re}}
\newcommand{\Imag}{\text{Im}}

\definecolor{darkgreen}{RGB}{0,128,0}
\newcommand{\review}[1]{#1}

\begin{document}

\maketitle

\begin{abstract}
A systematic mathematical framework for the study of numerical algorithms would allow comparisons, facilitate conjugacy arguments, as well as enable the discovery of improved, accelerated, data-driven algorithms.
Over the course of the last century, the Koopman operator has provided a mathematical framework for the study of dynamical systems, which facilitates conjugacy arguments and can provide efficient reduced descriptions. More recently, numerical approximations of the operator have enabled \review{ the analysis of a large number of deterministic and stochastic} dynamical systems in a completely data-driven, essentially equation-free pipeline.
\review{Discrete or continuous time numerical algorithms (integrators, nonlinear equation solvers, optimization algorithms) are themselves dynamical systems.
In this paper, we use this insight to leverage the Koopman operator framework} in the data-driven study of such algorithms and discuss benefits for analysis and acceleration of numerical computation.
For algorithms acting on high-dimensional spaces by quickly contracting them towards low-dimensional manifolds, we demonstrate how basis functions adapted to the data help to construct efficient reduced representations of the operator.
Our illustrative examples include the gradient descent and Nesterov optimization algorithms, as well as the Newton-Raphson algorithm.
\end{abstract}



\section{Introduction}

The relation between algorithms and dynamical systems has been studied for a long time. Due to the breadth of the topic, we cannot give an exhaustive review of the literature and will only briefly discuss some of the research areas important to this work. The book by Stuart and Humphries~\cite{stuart-1996} and the article by Chu~\cite{chu-2008} provide a broad overview of the relation between algorithms and dynamical systems.
\review{Much recent research} focuses on the discrete steps in the solution to linear systems and the relation to continuous dynamics on manifolds of matrices~\cite{chu-2008,guglielmi-2011,maeda-2016,miyatake-2018}.
In particular, algorithms for eigenvalue problems have been studied in this context, with the relation to integrable systems (Lax-pair formulation)~\cite{sun-2010}.
Studying the connection between algorithms and dynamical systems led to many surprising results, such as the definition of systems that can ``sort numbers''~\cite{brockett-1991}, or, more generally, solve optimization problems with constraints~\cite{gladkikh-2016}.
Chaotic dynamical systems that ``paint'', i.e., have stationary distributions according to the distribution of paint in an image, have been studied in relation to Markov Chain Monte Carlo algorithms~\cite{sahai-2017}.
Another connection relates particle dynamics to algorithmic solutions~\cite{edvardsson-2015}.
In the neural network research community, the idea of a (deep) neural network representing an iterated map of a discrete dynamical system has gained attention, particularly with respect to the continuous generator of this discrete map~\cite{pazos-2012}.
In general, optimizing the same loss function through different algorithms such as gradient descent or ``stochastic'' gradient descent may result in different trajectories, and different local or even global minima may be found.
Optimization algorithms modeling second-order dynamical systems (for example, Nesterov's method, sometimes called the ``heavy ball with friction'' algorithm) are favoured over traditional gradient descent methods, because of their ability to overcome local minima~\cite{bhaya-2009,noroozi-2009}.


The work of Koopman and von Neumann on ergodic, mixing, and chaotic dynamical systems~\cite{koopman-1932,neumann-1932} revealed that that there is a canonical, linear operator associated to each system. The operator acts on complex-valued functions of the system state, \review{the result being evaluation of the function at} a future time.
Linearity of the operator makes it amenable to finite-dimensional matrix approximations on computers, which was heavily exploited in the last twenty years~\cite{mezic-2005,budisic-2012}.
The Koopman operator and its adjoint, the Frobenius-Perron operator, are also known to provide optimal basis functions for uncertainty quantification~\cite{berry-2015d,giannakis-2017}.


In this paper, we  discuss benefits of the operator viewpoint on algorithms, when considering them as dynamical systems.
\review{The main problem we are addressing with the new operator viewpoint is the data-driven analysis of complex algorithms acting on high-dimensional state spaces. The Koopman operator and related numerical approximation methods for it allow us to treat this challenge in a unified way, solving many pressing issues such as acceleration and \review{formulation of data-driven} surrogate models, discovery of (almost) invariant sets and regions of convergence, and high-dimensional state spaces through model order reduction.
We demonstrate that a numerical approximation of the operator may be possible, even for systems with partially available or high-dimensional data.
We also show how properties of the Koopman operator, such as ergodic quotients and spectral analysis, offer valuable insight into algorithmic behavior.
Even if the only available data is a set of randomly distributed, single-iteration pairs of states, the operator \review{(and thus a sense of the algorithm)} can be approximately constructed.
We demonstrate a particularly novel approach to algorithm analysis through the computation of the spectrum of the Koopman operator. Summarizing a complex, nonlinear algorithm in such a meaningful way--and in the language of linear algebra--may provide a missing link to more fundamental results.
An example is the comparison of algorithms based on their Koopman spectrum, which may lead to a definition of metrics or distances between algorithms.
Conjugacy arguments have also been used to analyze dynamical systems for a long time. The operator viewpoint on algorithms can bring this successful approach to algorithm analysis, comparison, and development, even in a data-driven setting~\cite{bollt-2018b}.}

The remainder of the paper is organized as follows.
In section~\ref{sec:algorithms as dynamical systems}, we introduce the relation between algorithms and dynamical systems. In section~\ref{sec:koopman operator}, we describe the Koopman operator framework, followed by section~\ref{sec:edmd} with the Extended Dynamic Mode Decomposition algorithm for numerical approximation of the operator.
Acceleration and domain decomposition of algorithms is outlined in section~\ref{sec:acceleration}.
In section~\ref{sec:continuous nesterov}, we describe how to construct a continuous version of the iterative Nesterov algorithm in a data-driven way. 
Section~\ref{sec:high-dimensional state spaces} describes numerical approximations of the Koopman operator for algorithms that act on high-dimensional spaces, but move their state close to a low-dimensional manifold after a few iterations.
We analyze the Newton-Raphson method for root finding of polynomials on the complex plane in section~\ref{sec:newton method}, with explicit construction of the spectrum and eigenfunctions, as well as a numerical analysis of chaotic behavior of the method on the real line.
Appendix~\ref{sec:partial information} contains some preliminary computational experiments on the impact 
(on the numerical construction of the Koopman operator) of only partial information about the state space (due to finite sample size). 
We conclude with a summary of the results and an outlook on data-driven algorithm analysis in section~\ref{sec:conclusions}.

\section{The Koopman operator of algorithms}
The focus of this paper is on potential benefits of the data-driven study of algorithms as dynamical systems from the Koopman operator viewpoint.
In this section, we describe the necessary mathematical framework for dynamical systems, state the algorithms that we will use as examples, and briefly introduce the Koopman operator.
%
%


\subsection{Algorithms are dynamical systems}\label{sec:algorithms as dynamical systems}
In this paper, we consider algorithms that act on their state either continuously or in discrete steps (iterations).
We assume that the algorithms' state spaces $X\subseteq\mathbb{R}^d$ are smooth, $k$-dimensional Riemannian submanifolds embedded in Euclidean space of dimension $d\in\mathbb{N}$, $d\geq k$. The Riemannian metric on $X$ is induced by the embedding.
We define iterative algorithms as differentiable maps $a:X\to X$, with iteration number $n\in\mathbb{N}$. For a single iteration of the map $a$ on a state $x_n$, we write
\begin{equation*}\label{eq:discrete alg}
x_{n+1}=a(x_n),\ n\in\mathbb{N}.
\end{equation*}
For some algorithms, it is useful to consider continuously evolving state variables. In this case, the algorithm is typically represented as a differentiable vector field $v:X\to\mathbb{R}^k$. Starting from a given initial state $x_0\in X$, this vector field generates a curve $c:\mathbb{R}^+\to X$ in the state space through
\begin{equation}\label{eq:continuous alg}
\frac{d}{ds}c(s)=v(c(s)),\ c(0)=x_0.
\end{equation}
\Eqnref{continuous alg} describes the action of a continuous algorithm.
It can be reformulated as a discrete map through the definition of a flow: The map $\flow:\mathbb{R}^+\times X\to X$ is called a (semi-)flow of the vector field $v$ if, for all $t, s\in \mathbb{R}^+$ and all $x\in X$,
\begin{equation*}
\flow(t+s,x)=\flow(t, \flow(s, x)),\ \left.\frac{d}{dt}\flow(t,x)\right|_{t=0}=v(x).
\end{equation*}
If we fix a time interval $\Delta t\in\mathbb{R}^+$, the map $\flow(\Delta t,\cdot)=:\flow_{\Delta t}(\cdot):X\to X$ is an iterative algorithm (a discrete-time dynamical system). This construction is always possible: however, it is important to note that in general it is not possible to find a continuous algorithm $v$ for a given iterative algorithm $a$. We will discuss the necessary properties of $a$ for a continuous version to exist, and show how it can be approximated in a data-driven procedure through the consideration of the Koopman operator.
In this paper, we study the prototypical forms of algorithms listed in Table~\ref{tab:algorithms}, formulated in continuous and discrete time.

\begin{table}[htp]
\caption{\label{tab:algorithms}Algorithms formulated as dynamical systems, with their continuous (left) and discrete versions (right). Starting from a point $x_0\in\mathbb{R}^n$, the function $f:\mathbb{R}^n\to\mathbb{R}$ (or its absolute value, in case of root-finding algorithms) is to be minimized, where $\nabla f$ denotes its gradient,  and $H_xf$ the Hessian matrix at $x$. If $f$ is vector-valued, $J_xf$ denotes its Jacobian matrix at $x$. The symbol $h\in\mathbb{R}^+$ denotes a small, positive constant.}
\begin{tabular}{|p{0.45\textwidth}p{0.45\textwidth}|}\hline&\\
\multicolumn{2}{|c|}{Gradient Decent}\\&\\
$\dot{x}=-\nabla f(x)$&$x_{n+1}=x_n-h \nabla f(x_n)$\\&\\
\hline&\\
\multicolumn{2}{|c|}{Davidenko / Newton-Raphson method for optimization}\\&\\
$(H_xf)\dot{x}+\nabla f(x_0)=0$&$x_{n+1}=x_n-(H_{x_n}f)^{-1}\nabla f(x_n)$\\&\\
\hline&\\
\multicolumn{2}{|c|}{Davidenko / Newton-Raphson method for root finding}\\&\\
$(J_xf)\dot{x}+f(x_0)=0$&$x_{n+1}=x_n-(J_{x_n}f)^{-1} f(x_n)$\\&\\
\hline&\\
\multicolumn{2}{|c|}{Nesterov method}\\&\\
$\ddot{x}+\frac{r}{t}\dot{x}+\nabla f(x)=0, r\geq 3$&$x_{n+1}=y_n-h\nabla f(y_n)$\\
&$y_{n+1}=x_{n+1}+\frac{n}{n+3}(x_{n+1}-x_n)$\\
\hline
\end{tabular}
\end{table}

\subsection{Introduction to the Koopman operator}\label{sec:koopman operator}

\review{We start with a definition of the family of Koopman operators in the deterministic setting in section~\ref{sec:deterministic}. Details about the operator in the stochastic setting, and the relation with its adjoint, the Frobenius--Perron operator, are discussed in section~\ref{sec:stochastic setting}.}
\subsubsection{Deterministic setting}\label{sec:deterministic}
Given a deterministic dynamical system with vector field $v:X\to\mathbb{R}^k$, flow $\flow_t: X\to X$, initial condition $x_0\in X$, and 
\begin{equation}\label{eq:general vector field}
\left.\frac{d}{dt}\flow_t(x)\right|_{t=0}=v(x),\ S_0(x)=x_0,
\end{equation}
the family of {Koopman operators} $\K^t$ indexed by $t\in\mathbb{R}$ is a family of linear operators acting on the function space $\mathcal{F}$ of observables $\obs: X\to\mathbb{C}$ such that
\begin{equation*}
[\K^t\obs](x):=(\obs\circ \flow_t)(x).
\end{equation*}
The choice of function space $\mathcal{F}$ crucially determines the properties of $\K^t$. A typical choice is the space of complex-valued functions on the domain $X$ that are square-integrable with respect to a measure $\mu$, denoted
\begin{equation}\label{eq:l2 space}
\mathcal{F}=L^2(X,\mathbb{C},\mu):=\{\obs:X\to\mathbb{C}, \text{s.t.}\int_X |\obs(x)|^2d\mu(x)<\infty\}.
\end{equation}
The space $L^2(X,\mathbb{C},\mu)$ is typically equipped with an inner product
\begin{equation*}
\left\langle g_1,g_2\right\rangle:=\int_X g_1(x)\overline{g_2(x)}d\mu(x),
\end{equation*}
where the bar over $g_2(x)$ denotes complex conjugation.
The set of all $\K^t$ forms a continuous group w.r.t the strong operator topology. It is a semi-group if $\flow_t$ is not invertible for all $t$ (e.g., if the system approaches spatial infinity in finite time). The strong operator topology is defined through pointwise convergence on $(\mathcal{F},\|\cdot\|_\mathcal{F})$~\cite{engel-2006a}.
The map $t\mapsto \K^t$ is continuous in the strong operator topology on a compact subset $B\subset\mathbb{R}$ if, for $t\in B$,
\begin{equation*}
t\mapsto \K^tf
\end{equation*}
is continuous for every $f\in\mathcal{F}$. 
Continuous (semi-)groups are then defined as follows. 
The family $\{\K^t\}_{t\in\mathbb{R}}$ is a continuous (semi-)group of operators if, in the strong operator topology,
\begin{equation*}
\lim_{t\to 0}\|\K^t-I\|=0.
\end{equation*}
Typically, a fixed $t=t_0$ is chosen, and the associated Koopman operator $\K^{t_0}$ is studied. In the remaining text, if the specific choice of time is not important in the context, we will choose $t=1$ and drop the superscript $t$ to simplify notation.
An important feature of the operator is the spectrum $\sigma(\K)\subset\mathbb{C}$, where we consider the atomic (or pure-point) part $\sigma_{pp}$ and the absolutely continuous part $\sigma_{ac}$ (see \cite{budisic-2012,mezic-2017} for discussions of the singularly continuous part).
We can use the spectral theory for linear operators~\cite{neumann-1950,mezic-2005} to write
\begin{equation}\label{eq:koopman operator decomposition}
\K=\sum_k \underbrace{\Keval_k}_{\in\sigma_{pp}} P_{\Keval_k} + \int_{\sigma_{ac}}\Keval \text{d}E(\Keval),
\end{equation}
where $P_{\Keval}$ and $E(\Keval)$ are projection operators mapping functions to their projection in the eigenspaces associated to $\Keval$. \Eqnref{koopman operator decomposition} decomposes the operator $\K^t$ at $t=1$. Following functional calculus, the spectrum associated to other values of $t$ can be obtained by exponentiation of $\Keval$, i.e.  if $\Keval\in\sigma(\K)$ then $\Keval^t\in\sigma(\K^t)$.
For an eigenvalue $\Keval\in\sigma_{pp}(\K)$, there are functions $\phi_{\Keval}\in\mathcal{F}$ such that
\begin{equation*}
[\K^t\phi_\Keval](x)=(\phi_\Keval\circ \flow_t)(x)=\Keval^t \phi_\Keval(x),\ t\in\mathbb{R}^+.
\end{equation*}
These functions are called \textit{eigenfunctions} associated to the eigenvalue $\Keval$.
To predict the time evolution of a function $\obs\in\text{span}{\{\phi_{\Keval,k}\}}\subset\mathcal{F}$ such that $\obs=\sum_k c_k\phi_{\Keval,k},\ c_k\in\mathbb{C}$, we can thus write
\begin{equation}\label{eq:koopman predict}
\K^t\obs=\K^t\sum_k c_{k}\phi_{\Keval,k}=\sum_k c_{k}\K^t\phi_{\Keval,k}=\sum_k c_{k}\Keval^t\phi_{\Keval,k}.
\end{equation}
The coefficients $c_k$ are sometimes called ``Koopman modes'' for the observable $g$, since they are used to reconstruct it. For many dynamical systems, the span of the eigenfunctions is a rich subspace of $\mathcal{F}$, and many observables can be defined through their Koopman modes.

For continuous-time dynamical systems, another important object related to the Koopman operator group $\K^t$ is the infinitesimal generator $\A$, where
\begin{equation}\label{eq:koopman generator}
\A \obs:=\lim_{t\to 0}\frac{\K^t \obs-\obs}{t}
\end{equation}
for all $\obs\in D(\A)\subset\mathcal{F}$. For continuous semi-groups, the domain $D(\A)$ is dense in $\mathcal{F}$~\cite{engel-2006a}.
The infinitesimal generator has the same eigenfunctions as any $\K^t$. For fixed $t>0$ and eigenvalues $\Keval^t$ of $\K^t$, the eigenvalues $\Aeval$ of $A_\K$ satisfy $\exp(t\Aeval)=\Keval^t$. 
The generator is important because it induces the following relation on the eigenfunctions $\phi_\Keval$:
\begin{equation}\label{eq:koopman pde}
\left\langle\nabla\phi_\Keval(x),v(x)\right\rangle_{\mathbb{R}^n}=\Aeval \phi_\Keval(x)
\end{equation}
for all $x\in\mathbb{R}^n$ where the eigenfunctions and the vector field are defined (see \cite{bollt-2018b} for a discussion).
Note that equations~\eqref{eq:koopman generator} and \eqref{eq:koopman pde} reveal the relation between the generator $\A$ and the vector field $v$. If we interpret $v$ as a map from a point $x\in X$ to the coefficients $v_i$ for basis vectors $\partial_i$ of the tangent space $T_x X$, then we can write
\begin{equation*}
v(x)=\sum_i^d v_i(x) \partial_i.
\end{equation*}
With this reformulation, the vector field is an object that acts on functions $\obs:X\to\mathbb{R}$ from the left, i.e. $v\obs=\sum_i^d v_i \partial_i \obs$, which is exactly what \eqnref{koopman generator} describes.
Also, formally, 
\begin{equation*}
\K^t=\exp\left(t \A\right).
\end{equation*}
In this sense, the vector field $v$ acts as the generator of the Koopman operator {(semi)}\-group.
The generator $\A$ can be used to obtain an approximation of the vector field $v$ from only discrete-time samplings of the dynamical system. This was used to identify vector fields from data of dynamical systems~\cite{mauroy-2019}, and we will use it to identify vector fields for discrete algorithms in section~\ref{sec:continuous nesterov}.

The identification of vector fields using the Koopman operator is a global approach, taking into account all the data to obtain the vector field at every point. A different, more localized idea of identifying continuous vector fields on discrete data (e.g., \cite{rico-martinez-1992,rico-martinez-1995,he-2016}) is becoming increasingly popular in the machine learning literature.

\review{
\subsubsection{Stochastic setting}\label{sec:stochastic setting}
We provide a brief introduction to the operator definitions in the stochastic setting. The manuscript mostly focuses on deterministic algorithms, and thus we refer to  \cite{schuette-2016,crnjaric-zic-2019,klus-2018} for more theoretical details and applications of the stochastic case. Note that some numerical algorithms used to approximate the Koopman operator in the deterministic setting can be analyzed using a stochastic approximation of the system that accounts for numerical errors~\cite{williams-2015,crnjaric-zic-2019}.

Different from the deterministic systems based on equation~\eqref{eq:general vector field}, in this section, we start with a time-homogeneous stochastic process $\{x_t\}_{t\in\mathbb{R}^+_0}$ defined on the space $X$ with a probability measure $\mathbb{P}$.
%
In this setting, the formulation of the Koopman operator and its adjoint, the Frobenius-Perron operator, can be stated in terms of the transition density function~\cite{klus-2018}.
\begin{definition}
The transition density function $p:\mathbb{R}^+_0\times X\times X\to\mathbb{R}^+_0$ of a process $\{x_t\}_{t\in\mathbb{R}^+_0}$ is defined by
\begin{equation}
\mathbb{P}[x_{t+\tau}\in A|x_t=x]=\int_A p(\tau,x,y)dy
\end{equation}
for all measurable sets $A\subseteq X$.
\end{definition}
Note that for a deterministic system with flow $\phi$, the transition density function collapses to a Dirac delta centered on the point $\phi(t,x)$:
\begin{equation}
p(t,x,y)=\delta_{\phi(t,x)}(y).
\end{equation}
For $1\leq p\leq \infty$, the spaces $L^p(X,\mu)$ denote the $p$-integrable functions w.r.t. $\mu$ (see equation~\eqref{eq:l2 space} for $p=2$). The Koopman operator and its adjoint can now be defined as follows.
\begin{definition}
Let $q(t,\cdot)\in L^1(X,\mu)$ be a probability density and $f(t,\cdot)\in L^\infty(X,\mu)$ an observable of the system with flow $\phi$. For a given time step $\Delta t\in\mathbb{R}^+$,  the Koopman operator $K^{\Delta t}:L^\infty(X,\mu)\to L^\infty(X,\mu)$ is defined by
\begin{equation}\label{eq:def koopman stochastic}
K^{\Delta t}f(t,x)=\int_X p(\Delta t,x,y)f(t,y)d\mu(y)=\mathbb{E}\left[f(t,x_{t+\Delta t})|x_t=x\right],
\end{equation}
 and the Frobenius--Perron operator $P^{\Delta t}:L^1(X,\mu)\to L^1(X,\mu)$ is defined by
\begin{equation}\label{eq:def pf stochastic}
P^{\Delta t}q(t,x)=\int_X p(\Delta t,y,x)q(t,y)d\mu(y).
\end{equation}
\end{definition}
The Koopman and Frobenius--Perron operators are adjoint w.r.t the inner product on a function space $\mathcal{F}$ if for all $f,g\in\mathcal{F}$,
\begin{equation}
\left\langle K^{\Delta t}f,g\right\rangle_{\mathcal{F}}=\left\langle f,P^{\Delta t}g\right\rangle_{\mathcal{F}}.
\end{equation}

If for a given stochastic dynamical system the transition function satisfies the detailed balance condition
\begin{equation*}
p(t,x,y)=p(t,y,x)
\end{equation*}
for all $t\in\mathbb{R}^+$ and all states $x,y\in X$, the Koopman operator is self-adjoint (and hence equal to the Frobenius--Perron operator) on $L^1(X)\cap L^\infty(X)$, which is a direct consequence of the definitions of the operators (equations (\ref{eq:def koopman stochastic}--\ref{eq:def pf stochastic})).

The Frobenius--Perron and Koopman operators are the solution operators of the forward (Fokker--Planck) and backward Kolmogorov equations, respectively~\cite[section 11]{lasota-1994}. That is why they are also often referred to as \textit{forward} and \textit{backward} (transfer) operators.

General Smoluchowski equations of a $d$-dimensional system are given by
\begin{equation}\label{eq:general smoluchowski}
\text{d}\text{x}_t=-D\nabla V(\text{x}_t)\text{d}t+\sqrt{2dD}\text{d}\text{W}_t,
\end{equation}
with potential $V:X\to\mathbb{R}$, diffusion coefficient $D$, and Wiener process $W_t$. They are time-reversible~\cite{mazo-2008,leimkuhler-2005}, and so their transition density function satisfies the detailed balance condition. Therefore, the two transfer operators are identical (self-adjoint). The fact that this holds for equations of type \eqref{eq:general smoluchowski} is extremely strong and useful, particularly for the analysis of molecular dynamics simulators. As many algorithms can be reformulated in this setting, too, these strong results can be carried over directly.
}


\subsection{Illustrative example: Koopman operator of the Euler method}
We use a typical example from the stability analysis of algorithms to familiarize the reader with the concept of Koopman operators of algorithms.
Consider the following ordinary differential equation (ODE) with parameter $a>0$,
\begin{equation}\label{eq:simple ode}
\left.\frac{d}{dt}\flow_t(x)\right|_{t=0} = -a x,\ x\in\mathbb{R},\ t\geq 0,
\end{equation}
and solution $\flow_t(x)=\exp(-a t)x$. A trajectory $\flow_t(x_0)$ can be approximated with a numerical algorithm, \review{an initial value solver,} given only the ODE~(\ref{eq:simple ode}) and an initial condition $x_0\in\mathbb{R}$. Here, we choose the forward Euler method, compute its Koopman operator, and show how stability analysis of the method is related to the spectrum of the operator.
For a given step size $\Delta t>0$, the forward Euler method is given through the following iterative scheme, starting at $x_0$ for $n=0$:
\begin{equation}
x_{n+1}=x_n + \Delta t \left.\frac{d}{dt}\flow_t(x_n)\right|_{t=0}.
\end{equation}
For the ODE~(\ref{eq:simple ode}), we thus have the linear, discrete-time dynamical system
\begin{equation}\label{eq:euler method}
x_{n+1}=x_n - a\Delta t x_n = (1-a\Delta t) x_n.
\end{equation}
We now choose the function space $\mathcal{F}=\text{span}\{x^k, x\in\mathbb{R}, k\in\mathbb{N}\}$, with the identity function $g(x)=x$ as the generator of the basis. The Koopman operator $\K^n:\mathcal{F}\to\mathcal{F}$ applied to the generator is
\begin{equation*}
[\K^n\obs](x)=\obs((1-a\Delta t)^nx) = (1-a\Delta t)^n\obs(x),
\end{equation*}
which shows that $\obs$ is an eigenfunction of $\K^n$ associated to the eigenvalue $(1-a\Delta t)$.
On the space $\mathcal{F}$, the spectrum of $\K$ consists of the eigenvalues $\omega_k= (1-a\Delta t)^k$, $k\in\mathbb{N}$, associated to the eigenfunctions $\phi_k(x)=\obs(x)^k=x^k$. It is easy to see that for $a\Delta t\in[0,2]$, $\K^n$ is a contraction operator on $\mathcal{F}$: for all $k$, $|\Keval_k|\leq 1$. This is a sufficient condition for the numerical stability of the Euler method in this region.
Note that system~(\ref{eq:simple ode}) has eigenvalues $\Keval_k=\exp(-ka)$ on $\mathcal{F}$, and so the discrete system (\ref{eq:euler method}) defined by the Euler method is not conjugate to the continuous system for any $\Delta t>0$ (for conjugacy, the eigenvalues would have to be the same). Consistency--and thus, convergence--of the algorithm can now also be interpreted in the spectral sense, by setting $n:=\lceil{1/\Delta t}\rceil$ (rounding to the next largest integer) and applying the formula
\begin{equation*}
\lim_{n\to\infty}\left(1-\frac{a}{n}\right)^{n}=\exp(-a).
\end{equation*}
After this pedagogical example, we show how Koopman operators can provide a useful framework to understand more general algorithms.
In a more applied setting, ``the algorithm'' is usually a large, complex piece of software that is treated as a black box data generator (for example, an optimization procedure for a supply chain). Then, numerical approximations of the Koopman operator can be used. In the next section, we briefly describe Extended Dynamic Mode Decomposition as an example for such an approximation procedure from data.

\section{Data-driven approximation of the Koopman operator}\label{sec:edmd}
We briefly cover the approximation methods employed in the computational experiments of later sections. The methods include an approximation of the Koopman operator with pure-point spectrum (Extended Dynamic Mode Decomposition), extracting the infinitesimal generator of the family of Koopman operators, as well as an approximation algorithm in case the Koopman operator is unitary and has continuous spectrum.

\subsection{Extended Dynamic Mode Decomposition}

One extension of Dynamic Mode Decomposition \review{(DMD)}~\cite{schmid-2010} to nonlinear systems with pure-point spectrum is Extended DMD~\cite{williams-2015,williams-2015b,williams-2015c} (EDMD).
Many more numerical methods exist, such as Generalized Laplace Averaging~\cite{mezic-2013,mauroy-2013}, EDMD with dictionary learning through neural networks~\cite{li-2017}, dictionary selection through $L^1$-optimization~\cite{brunton-2016c}, etc.

We consider the state space of a given system--usually, a smooth manifold $M$ embedded in Euclidean space $\mathbb{R}^d$--sampled with a finite number of points $X=\{x_k\in M\subset\mathbb{R}^d|k=1,\dots,N_X\}$, $N_X\in\mathbb{N}$.
The main idea of (Extended) Dynamic Mode Decomposition is to approximate the action of the operator  on a function space $\mathcal{F}$ over $M$ by choosing a finite-dimensional subspace spanned by a finite set of real-valued functions of $M$ (called a ``dictionary'', with ``observables'' as elements). For $\mathcal{F}$, we typically use the space $L^2(M,\mathbb{C},\mu)$, see \eqnref{l2 space}.
The EDMD algorithm computes the action of the Koopman operator on the dictionary at the points in $X$, and then approximates the Koopman operator by solving a least-squares problem:
\begin{enumerate}
\item Given a dictionary $D=\{d_k:M\to\mathbb{R}|k=1,\dots,N_D\}\subset \mathcal{F}$ with $N_D$ observables, construct
\begin{equation*}
G=D(X)=\left[d_1(X),d_2(X),\dots,d_{N_D}(X)\right]\in\mathbb{R}^{N_X\times N_D}.
\end{equation*}
\item For a fixed $t>0$, compute the action of the Koopman operator $\K^t$ on the dictionary elements using the flow map $\flow_t:M\to M$ of the given system:
\begin{equation*}
[\K^t d_{k}](x)=(d_k\circ \flow_t)(x),
\end{equation*}
and construct the matrix
\begin{equation*}
A=[\K^tD](X)=\left[\mathcal{{K}}^td_1(X),\K^td_2(X),\dots,\K^td_{N_D}(X)\right]\in\mathbb{R}^{N_X\times N_D}.
\end{equation*}
\item Approximate the operator $\K^t$ through the matrix
\begin{equation*}
K=\frac{1}{N_D^2}(G^T G)^\dagger (A^T A)\in\mathbb{R}^{N_D\times N_D},
\end{equation*}
\noindent where $(G^T G)^\dagger$ denotes the pseudo-inverse of the matrix $(G^T G)$.
\end{enumerate}
The choice of the dictionary $D$ is crucial for a successful approximation of the operator, see~\cite{li-2017,williams-2015,williams-2015c} and references therein. In this paper, we employ thin-plate radial basis functions~\cite{duchon-1977},
\begin{equation}\label{eq:thin plate}
d_x(y)=\|x-y\|^2\ln(\|x-y\|+\delta),\ x,y\in X,\ \delta>0,
\end{equation}
where the point $x$ is the center of the radial basis function and $\delta$ is a small, positive constant to extend the function to points $y=x$ (here: $\delta=10^{-3}$).
The centers $x$ are uniformly distributed over the data set $X$, where the uniform distribution is constructed through \review{an appropriate} k-means algorithm, clustering the data in $N_D$ clusters and using the centers of the clusters as the centers for the basis functions.

\subsection{Approximation of the infinitesimal generator}

After some of the eigenfunctions $\hat{\Kefunc}$ and eigenvalues $\Keval$ of $\K$ are available through their numerical approximations, we can represent the approximation of the infinitesimal generator $\A$ for the group $\K^t$ as a matrix $L=V\frac{\ln[\Lambda]}{\Delta t}V^{-1}$ (see~\cite{mauroy-2019}), where $V$ contain the eigenvectors of the matrix representation of a particular $\K^{\Delta t}$, obtained with EDMD from snapshot data $(x,\flow_{\Delta t}(x))$. For the coordinate functions $x_i\in\mathcal{F}$ and a finite step size $\Delta t$,
\begin{equation}\label{eq:Koopman generator action}
\A x_i = \lim_{t\to 0}\frac{(\K^t-I)x_i}{t} \approx  \frac{\ln\left[\Lambda\right]}{\Delta t} \hat{\phi}(x)\, C_i,
\end{equation}
where $\Lambda$ is the diagonal matrix of eigenvalues of $\K$, $\ln$ is the complex logarithm (where we pick the principal branch), and $C$ are the ``Koopman modes'' associated to the coordinate functions $x_i:X\to\mathbb{R}$. The approximation of \review{the vector field $v$ generating the flow $S_t$} is thus
\begin{equation}\label{eq:vector field from generator}
v(x) \approx  \frac{\ln\left[\Lambda\right]}{\Delta t} \hat{\phi}(x) C.
\end{equation}
Note that the approximation of the generator through the logarithm of the eigenvalues involves several challenges: numerical issues with eigenvalues of $\Lambda$ close to zero, and non-invertibility (or rather, non-uniqueness) of the complex logarithm. 
The issue with non-uniqueness is noted by Mauroy and Goncalves~\cite{mauroy-2019}. In the example where we construct $L$ (section~\ref{sec:continuous nesterov}), the spectrum of the Koopman operator lies in the unit disk, so the complex logarithm is unique when picking its principal branch. To mitigate numerical issues, we set all eigenvalues with $\Real{(\ln\Keval)}<-2/\Delta t$ to zero.

\subsection{Approximation of the continuous spectrum}\label{sec:approx continuous spectrum}

The EDMD algorithm is applicable for approximating the Koopman operator if its spectrum only consists of eigenvalues (pure-point spectrum). Not many numerical approximations are available \review{for operators with continuous spectra}. Here, we briefly describe the results from Korda, Putinar and Mezi\'{c}~\cite{korda-2018} for unitary Koopman operators, where the spectrum is concentrated on the unit circle $\mathbb{T}$ in the complex plane.
We will use this approximation to analyze the Newton algorithm in section~\ref{sec:newton method}.
If the Koopman operator of a system is unitary, we can write it as an integral over a projection operator-valued measure $E$ on the unit circle $\mathbb{T}\subset\mathbb{C}$,
\begin{equation*}
\K =\int_\mathbb{T}z dE(z).
\end{equation*}
For any observation function ${\obs}$ in the domain $\mathcal{F}$ of $\K$, the measure $E$ defines a real-valued, positive measure $\mu_{\obs}$ on $\mathbb{T}$ through
\begin{equation*}
\mu_{\obs}(A):=\left\langle E(A)\obs, \obs\right\rangle_\mathcal{F}
\end{equation*}
for all Borel sets $A\subset\mathbb{T}$. The moments $m_k$, $k\in\mathbb{Z}$ of $\mu_{\obs}$ are defined by
\begin{equation*}
m_k:=\int_\mathbb{T}z^kd\mu_{\obs}(z)=\left\langle \K^k \obs,\obs\right\rangle_\mathcal{F},
\end{equation*}
where the last identity follows from the spectral theorem.
The measure $\mu_{\obs}$ depends on the choice of observable $\obs$, but under certain conditions on $\obs$, $\mu_{\obs}$ fully determines the operator $\K$.
The conditions are satisfied if $\obs$ is $\star$-cyclic, meaning that repeated applications of $\K$ to ${\obs}$ span the function space $\mathcal{F}$, see~\cite{korda-2018}. In this case, approximating the measure $\mu_{\obs}$ not only reveals information about the specific observable $\obs$, but also about the underlying system (and its operator $\K$).
The density of $\mu_{\obs}$ on the unit circle can be visualized after a numerical approximation of the moments.
If we have $N$ observations $y_i=\obs(x_i)$ along a trajectory of an ergodic system, we can estimate the moments through
\begin{equation*}
m_k\approx\frac{1}{N-k}\sum_{i=1}^{N-k} y_{i+k} \bar{y}_{i}.
\end{equation*}
Define $\psi_N(z)=[1,z,z^2,\dots,z^N]^T$, and write the Christoffel–Darboux kernel as
\begin{equation*}
K_N(z,w)=\psi_N(w)^HM^{-1}\psi_N(z),
\end{equation*}
with a semi-positive definite, Hermitian Toeplitz matrix ${M}$ defined through
\begin{equation*}
M:=\left[
\begin{matrix}
m_0&\overline{m}_1&\cdots&\cdots&\cdots&\overline{m}_N\\
m_1&m_0&\overline{m}_1&\ddots&\ddots&\overline{m}_{N-1}\\
\vdots&m_1&{m_0}&\ddots&\ddots&\overline{m}_{N-2}\\
\vdots& & \ddots& \ddots&\ddots &\vdots\\
\vdots&&&\ddots&m_0&\overline{m}_{1}\\
m_N&\cdots&\cdots&\cdots&m_1&{m}_{0}\\
\end{matrix}
\right].
\end{equation*}
Inverting the matrix $M$ is non-trivial, see~\cite{korda-2018} for details. Instead of inverting $M$ directly, the authors propose to invert a matrix $\tilde{M}$ where the elements $m_k$ in $M$ are replaced by elements $$\tilde{m}_k:=\left\lbrace\begin{matrix}m_k+1&\text{if }k=0,\\m_k&\text{if }k> 0.\end{matrix}\right.$$
Given the kernel $\tilde{K}_N$ (now using $\tilde{M}$ instead of $M$), the continuous density $\rho$ of $\mu_{\obs}$ at a given angle $\theta\in\mathbb{T}\cong [0,2\pi)$ can be approximated through
\begin{equation*}
\frac{N+1}{\tilde{K}_N(\exp(i\theta),\exp(i\theta))}-1\approx\rho(\theta).
\end{equation*}



\section{Analysis of algorithms through their Koopman operator}\label{sec:examples}

This section contains illustrative examples of algorithm analysis through their Koopman operator. 
\review{We demonstrate several aspects of our solution to the main challenge: how to analyze complex, nonlinear algorithms in a data-driven setting within a unified framework. In particular, we show
\begin{enumerate}
    \item how to accelerate algorithms by constructing data-driven surrogate models,
    \item approximations of infinitesimal generators for discrete-time algorithms,
    \item approximations of the (asymptotic, long-time) operator in high-dimensional state spaces,
    \item and a combination of an explicit operator analysis for the Newton-Raphson method with a data-driven approximation of its continuous spectrum.
\end{enumerate}
}

\subsection{Acceleration}\label{sec:acceleration}
We can employ the predictive capabilities of the Koopman operator framework to accelerate the application of algorithms. First, the Koopman operator eigenvalues, eigenfunctions, and modes are constructed in an offline phase, using EDMD (see \secref{edmd}). Second, given a new initial condition close to the ones in the data set used for constructing the operator, we can use the approximation as a data-driven surrogate model for the original algorithm.

Such a surrogate model may be useful in accelerating the training of neural networks. There, the state space is often very high-dimensional, because each state consists of all weights and biases of the network. The typical optimization algorithm, stochastic gradient descent, employed to change the state in order to minimize a certain loss function, is based on the gradient descent algorithm we study in this section. \review{The stochastic version can be analyzed similarly, using the setting described in section~\ref{sec:stochastic setting}.} Different types of training (e.g. through mini-batches, or dropout) will induce different biases, and the Koopman operator based surrogate may provide a path towards quantifying \review{and comparing such inductive biases}.

\begin{figure}[htp]
\centering
\begin{tabular}{cc}
\includegraphics[width=\figureThirdWidth]{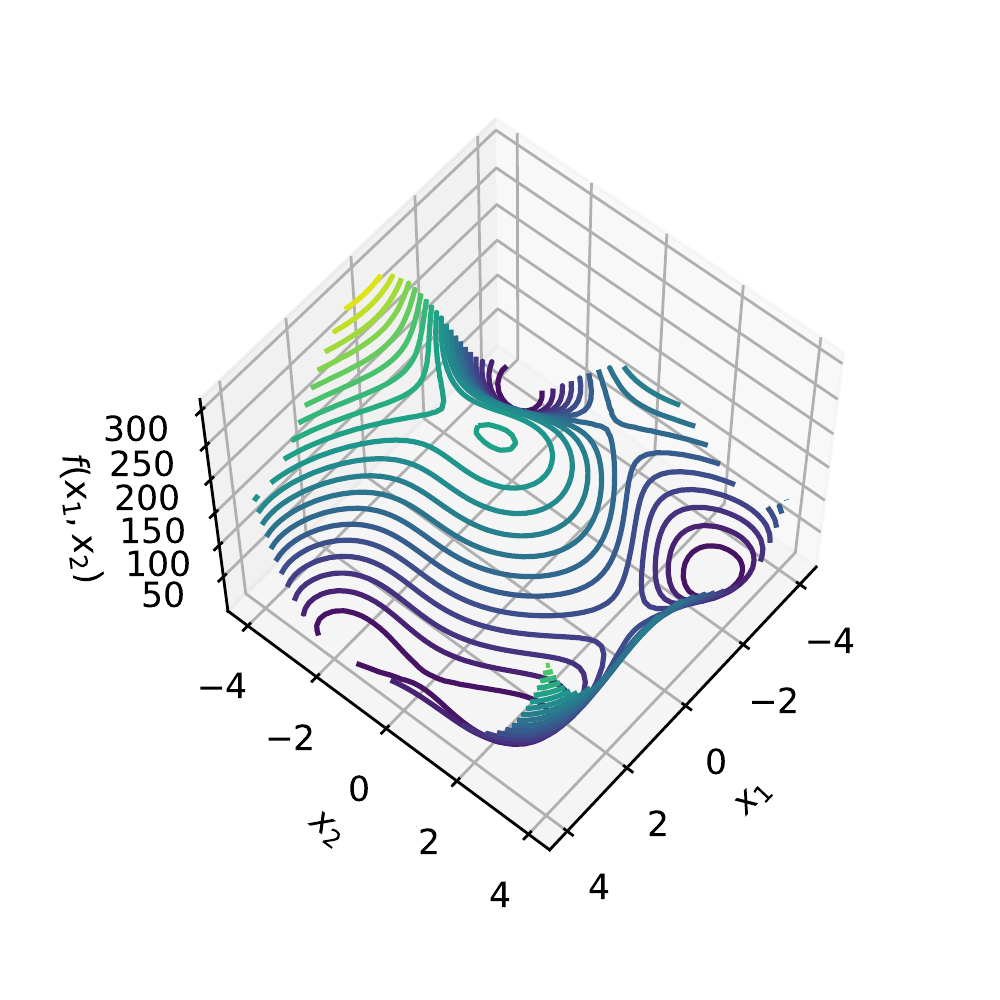}&
\includegraphics[width=\figureThirdWidth]{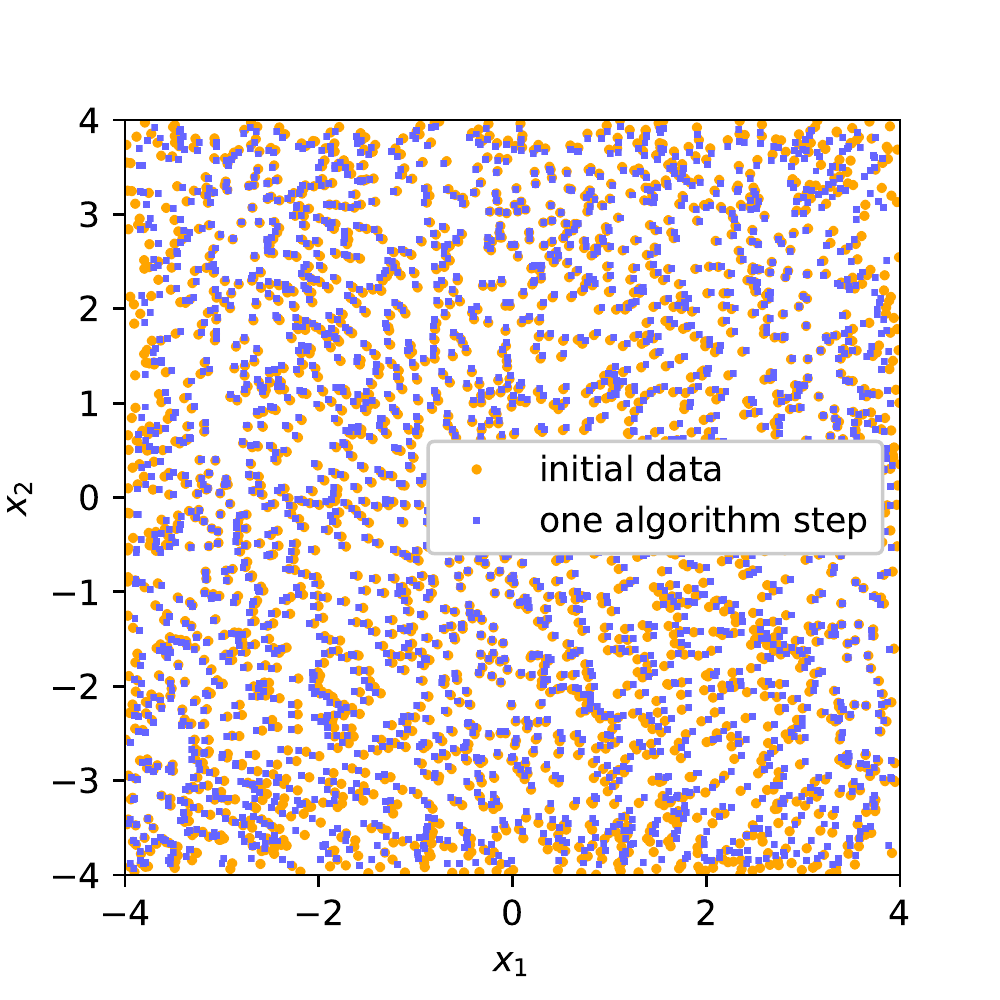}\\
(a)&(b)
\end{tabular}
\caption{\label{fig:himmelblau himmelblau_data_on_phase_space}Level sets of Himmelblau's function (a), and sample data $(y_n,y_{n+1})$ from Gradient Descent (b), where $y_{n+1}$ is generated through one iteration step with $\Delta t=.001$.}
\end{figure}

\begin{figure}[htp]
\centering
\begin{tabular}{ccc}
\raisebox{.35cm}{\includegraphics[width=\figureThirdWidth]{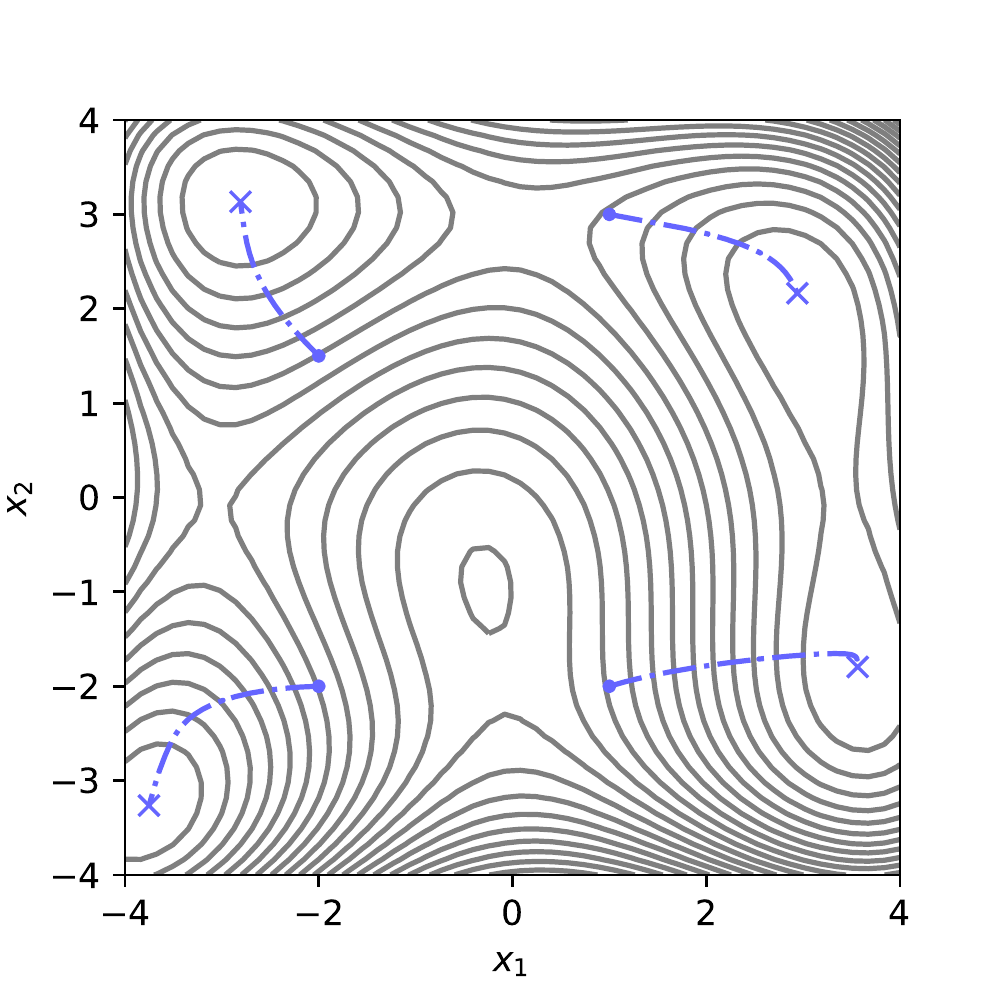}}&
\includegraphics[width=\figureThirdWidth]{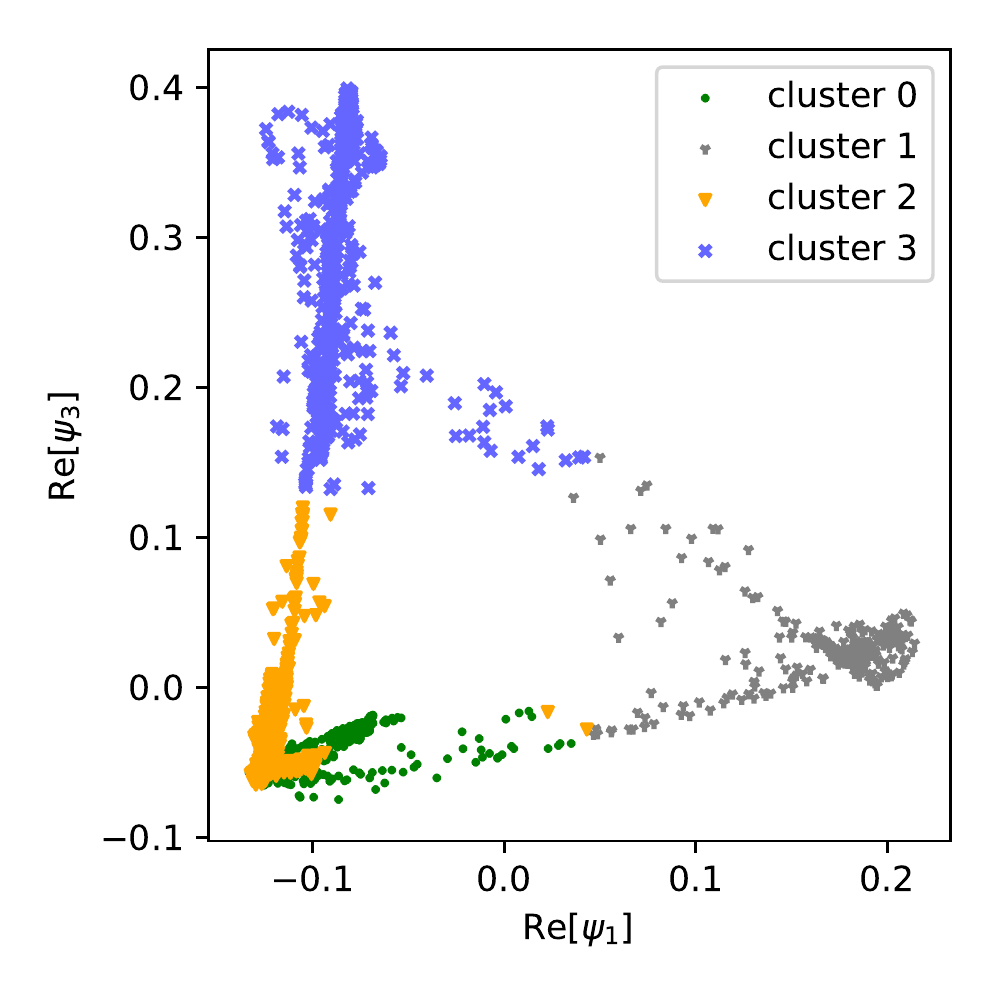}&
\includegraphics[width=\figureThirdWidth]{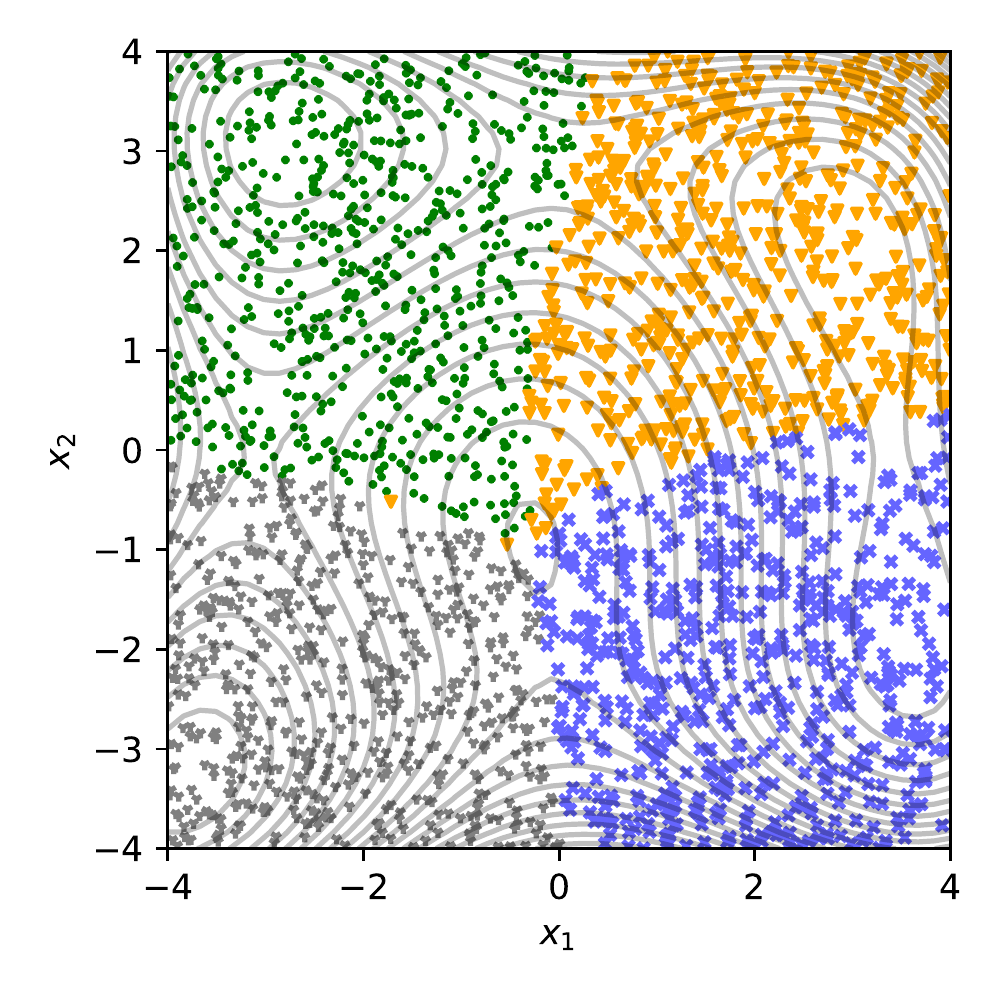}\\
(a)&(b)&(c)
\end{tabular}
\caption{\label{fig:himmelblau himmelblau_trajectories_on_phase_space} Panel (a) shows predicted trajectories from four initial points. All four trajectories approach the correct minima in their respective basin of attraction. An ergodic decomposition of the state space through k-means clustering of values from four complex-valued Koopman eigenfunctions associated to eigenvalues close to 1 (a projection of the eight-dimensional data to two dimensions is shown in panel b) separates the four basins of attraction (c).}
\end{figure}

\Figref{himmelblau himmelblau_trajectories_on_phase_space} shows the setting used for the example in this section, where we study the Gradient Descent algorithm minimizing Himmelblau's function~\cite{himmelblau-1972}---a standard test function for optimization algorithms:
\begin{equation*}
f(x_1,x_2)=(x_1^2+x_2-11)^2+(x_1+x_2^2-7)^2.
\end{equation*}
This function has one local maximum and four local minima with three connecting saddle points (see \figref{himmelblau himmelblau_data_on_phase_space}, a).
As a dictionary for EDMD, we employ 500 thin-plate radial basis functions (see \eqnref{thin plate}), chosen uniformly distributed on $M=[-4,4]^2$.
In addition to the radial basis functions, we add the two coordinate functions $x_1,x_2$ as well as the constant function with value $1$ to the dictionary.

\Figref{himmelblau himmelblau_trajectories_on_phase_space} shows four predicted trajectories using the approximated Koopman operator.
To generate a trajectory, we follow the prediction steps outlined in \secref{koopman operator} (\eqnref{koopman predict}), using all 503 eigenvalues and eigenfunctions.

The eigenfunctions of $\K^t$ at eigenvalue $\Keval=1$ are special, as they are preserved under the flow of the system---they do not decay, expand, or oscillate. This can be used to construct an \textit{ergodic decomposition}~\cite{budisic-2009,budisic-2012} of the state space, separating the basins of attraction. In the system discussed here, we expect to see four basins of attraction, corresponding to four attracting fixed points, and accordingly, at least four eigenfunctions associated to eigenvalue $\Keval=1$. To approximate the ergodic decomposition, we consider the map $M\ni p\mapsto [\Kefunc_1(p),\dots,\Kefunc_4(p)]^T\in\mathbb{C}^4$, i.e. from the states $M$ into the values of eigenfunctions associated to $\Keval=1$. In this space (shown projected to two dimensions in \figref{himmelblau himmelblau_trajectories_on_phase_space}, b), there should be four clusters of values corresponding to the different values the eigenfunctions take on over the corresponding basins of attraction. Clustering the data using the k-means algorithm then allows us to distinguish the original points in $M$ (plot c in \figref{himmelblau himmelblau_trajectories_on_phase_space}).

\subsection{Continuous versions of discrete algorithms}\label{sec:continuous nesterov}
\begin{figure}[htp]
\centering
\includegraphics[width=\figureFullWidth]{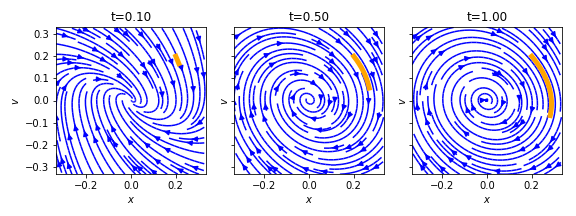}
\includegraphics[width=\figureFullWidth]{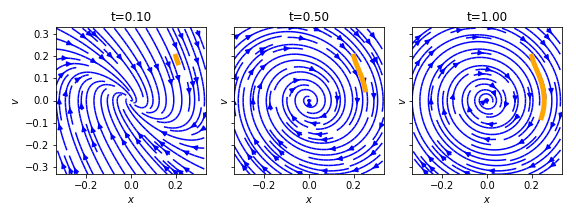}
\caption{\label{fig:nesterov edmd}Approximation of the vector field through EDMD (top row), and the ODE found by Cand\`{e}s et al.~\cite{su-2014} (bottom row), at times $t\in \{\frac{1}{10}, \frac{1}{2}, 1\}$. The orange trajectory shows an example solution over this (time-dependent!) vector field.}
\end{figure}

In recent years, extensive research has focused on re-deriving established algorithms as discrete versions of continuous dynamical systems, in this way ``explaining'' the original, iterative algorithm.
In the case of the Nesterov method, the appropriate continuous system derived by  Cand\`{e}s et al.~\cite{su-2014}  is a time-dependent, second-order ordinary differential equation. This ODE can be (trivially) transformed to a first-order system of autonomous equations, given through
\begin{equation}\label{eq:nesterov ode a 1}
\dot{x}=v,\ 
\dot{v}=-\frac{r}{t}v-\nabla f(x),\ 
\dot{t}=1,
\end{equation}
where $\dot{x}$ now refers to a derivative with respect to a new time variable. 
In the example we discuss below, the function to minimize is $f(x)=\frac{1}{2}x^2$, $x\in \mathbb{R}$. As in the previous section, we use EDMD to approximate the eigensystem of the Koopman operator with 125 thin-plate radial basis functions (see \eqnref{thin plate}) as well as the $x,v,t$ coordinates and the constant function as a basis for the function space of observables on $\mathbb{R}^3$ ($(x,v,t)$-space).
The operator approximation is performed with data $(z_n,z_{n+1})$ from the discrete algorithm with step size $h=0.01$ (see table~\ref{tab:algorithms}), where $z_n:=(x(n),(x(n)-x(n-1))/h,t(n))$.
We then approximate the corresponding vector field on the three coordinates of $z$ through equations \eqref{eq:Koopman generator action} and \eqref{eq:vector field from generator}, which should be similar to the vector field in \eqnref{nesterov ode a 1}.
\Figref{nesterov edmd} shows streamlines of the vector field~\eqref{eq:nesterov ode a 1} for three different values of $t$, with a distinguished sample trajectory, compared to the vector field approximated by using the Koopman operator.
This example illustrates that it is possible to extract a continuous version (an infinitesimal generator) of a discrete algorithm.

\subsection{High-dimensional state spaces}\label{sec:high-dimensional state spaces}

Many algorithms are operating on points in high-dimensional state spaces that are  pa\-ra\-me\-trized  through a large number of variables. In this case, sampling the full state space to obtain an accurate approximation of the Koopman operator for the algorithm is difficult, especially when the action of the algorithm on the state is nonlinear. However, if a few iterations of the algorithm quickly contract the whole state space towards a low-dimensional subset, then numerical approximations can be successful \review{in reducing the overall computational effort}.
The EDMD algorithm does not need an explicit expression for this low-dimensional subset if, for example, radial basis functions are used for the dictionary.

\begin{figure}
\centering
\includegraphics[width=\figureHalfWidth]{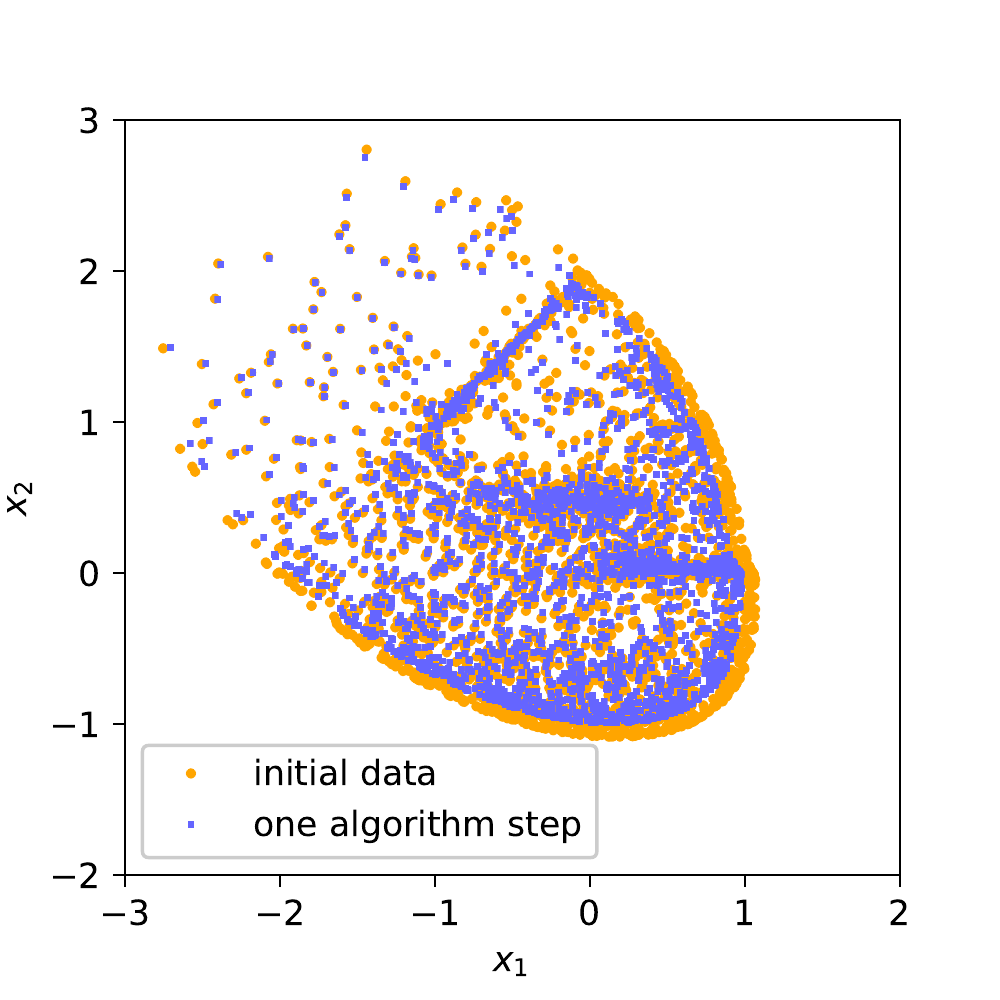}
\raisebox{-.06cm}{\includegraphics[width=4.25cm]{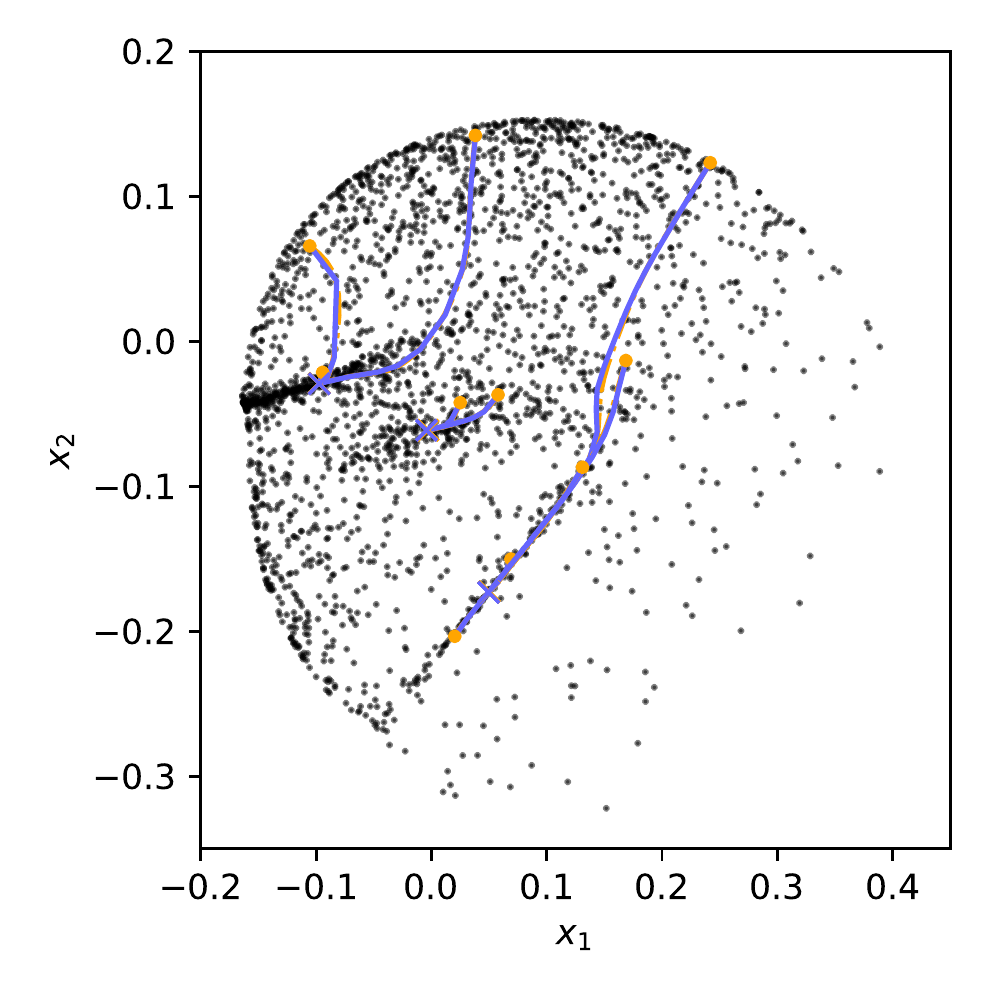}}
\includegraphics[width=\figureFullWidth]{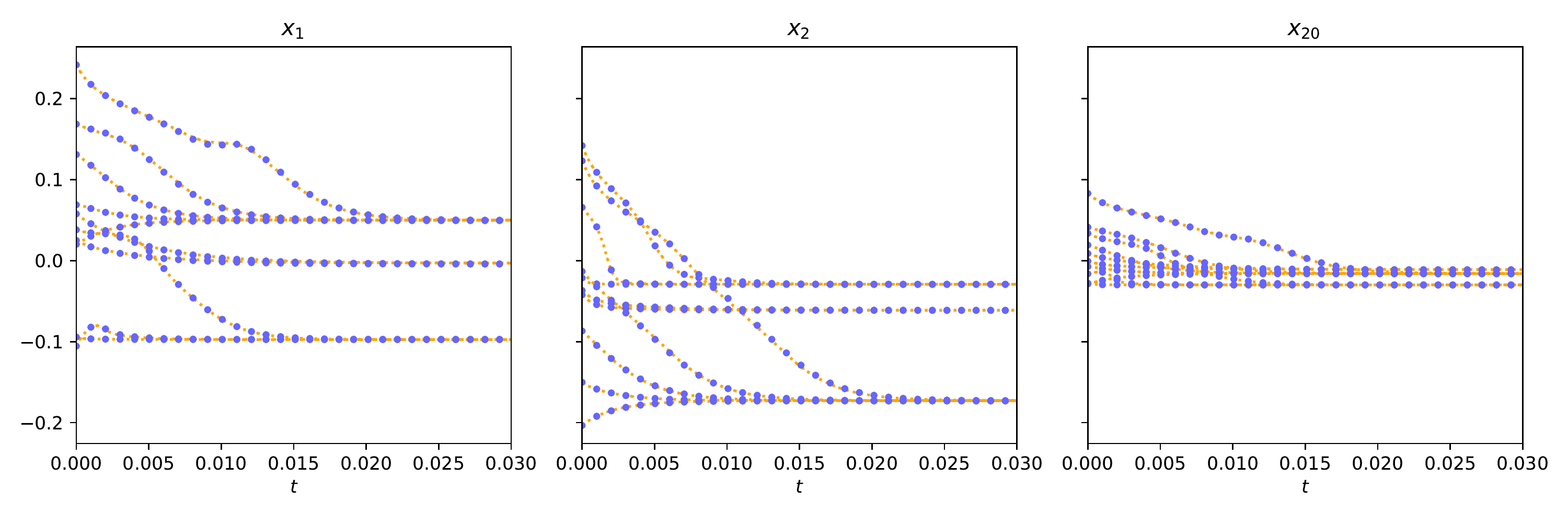}
\caption{\label{fig:mb1}Top-Left: Data gathered by sampling the 100-dimensional state space with a Gaussian distribution, converging to the low-dimensional manifold in five iterations of gradient descent (end points marked in orange), and then recording one additional iteration (blue points). The coordinates $x_1$ and $x_2$ here are only the first two of all 100 present in the data set.
Top-Right: Predicted trajectories of the algorithm (blue) that were computed after the construction of the operator, for ten initial conditions chosen at random from the initial data set (orange points). The bottom panels show predicted (blue, dots) and actual (orange, dashed) trajectories for the coordinates $x_1$, $x_2$, and $x_{20}$, to demonstrate that all coordinate trajectories are reconstructed accurately.}
\end{figure}

As an example, we construct a high-dimensional optimization problem with a low-\-di\-men\-sio\-nal, attracting manifold, where gradient descent is used as an optimizer. On the manifold, starting from the M\"uller-Brown potential, we add quadratic terms to obtain asymptotically stable, attractive behavior in the directions transverse to the manifold. On the manifold (here, a two-dimensional subspace of $\mathbb{R}^{100}$ parametrized by $x,y$), the potential is defined by
\begin{equation*}
V(x,y)=\sum_{k=1}^4 A_k \exp\left( a_k(x-x_k^0)^2+b_k(x-x_k^0)(y-y_k^0)+c_k(y-y_k^0)^2 \right),
\end{equation*}
with
\begin{eqnarray*}
A&=&(-200,-100,-170,15),\\
a&=&(-1,-1,-6.5,0.7),\\
b&=&(0,0,11,0.6),\\
c&=&(-10,-10,-6.5,0.7),\\
x^0&=&(1,0,-0.5,-1),\\
y^0&=&(0,0.5,1.5,1).
\end{eqnarray*}
The cost function $C:\mathbb{R}^{100}\to \mathbb{R}$ on the embedding space is then defined through
\begin{equation*}
C(\textbf{x})=V([U\textbf{x}]_1, [U\textbf{x}]_2) + \sum_{i=0}^{98} (U\textbf{x})_i^2, \ \textbf{x}\in\mathbb{R}^{100},
\end{equation*}
where we rotate the data by a random, unitary matrix $U\in\mathbb{R}^{100\times 100}$ into the 100 dimensional space to demonstrate that the EDMD algorithm with radial basis functions does not depend on the chosen embedding coordinate system.
Now, all coordinates have non-trivial dynamics, but the low-dimensional, attractive manifold is still present.
To approximate the operator, we use 625 thin-plate radial basis functions with centers randomly distributed over the initial data set in the 100-dimensional space, as well as the 100 coordinates $x_1,\dots,x_{100}$ and the constant function (i.e., the final dictionary has 726 elements).
\review{The initial data consists of 2500 data points, sampled in a 100-dimensional standard normal distribution around zero, and then evolved forward with five iterations of the algorithm to converge towards the low-dimensional structure. \textit{Note that this number of points is--by far--not enough to densely sample the high-dimensional space,} but the attracting behavior of the algorithm towards the low-dimensional structure is nonetheless sufficient to numerically approximate the operator. We use radial basis functions for EDMD to circumvent the need for an explicit approximation of the low-dimensional space: the functions only depend on the distance between the points, not the ambient dimension.}
As in \secref{acceleration}, we can again accurately predict trajectories of gradient descent (see \figref{mb1}).

\subsection{Koopman operators for the Newton-Raphson method}\label{sec:newton method}
\begin{figure}[ht!]
\centering
\includegraphics[width=.35\textwidth]{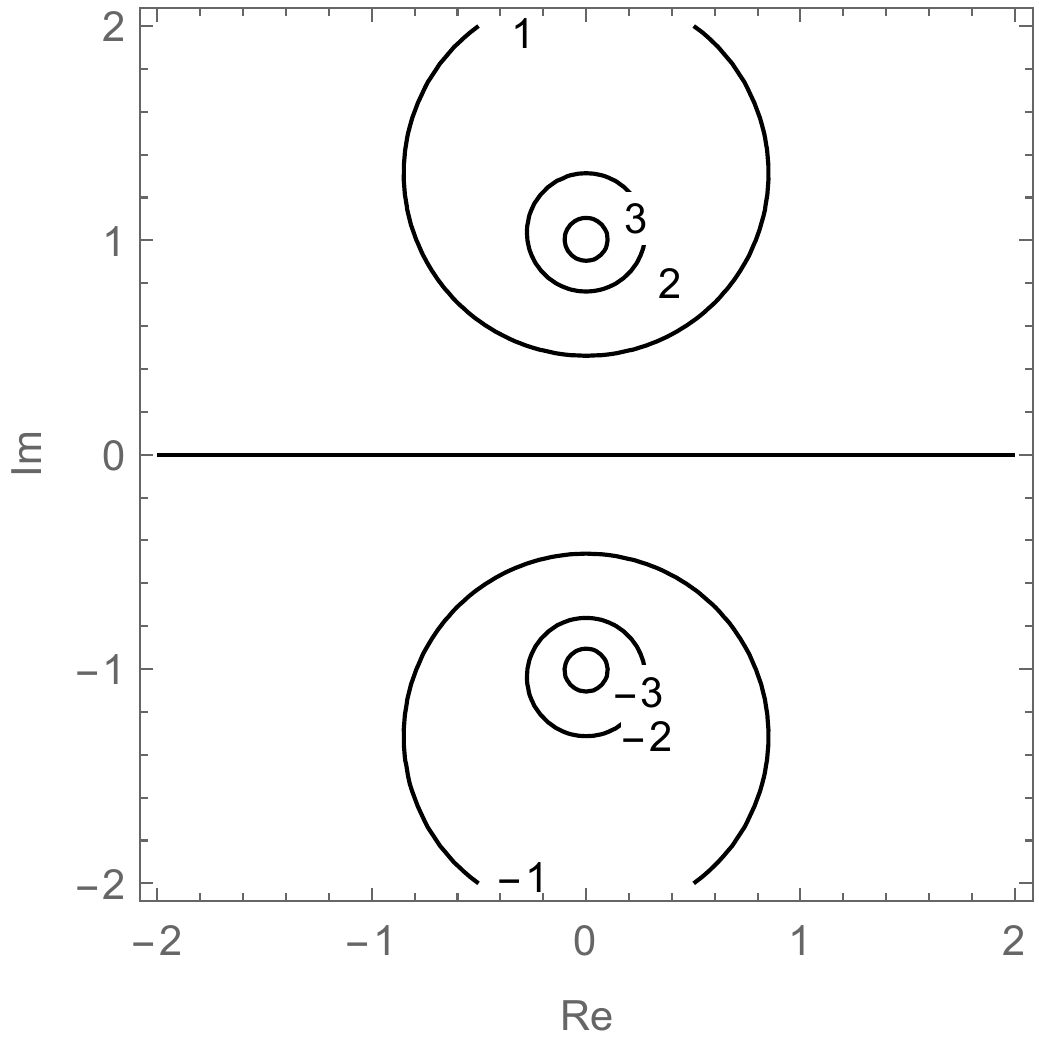}
\includegraphics[width=.55\textwidth]{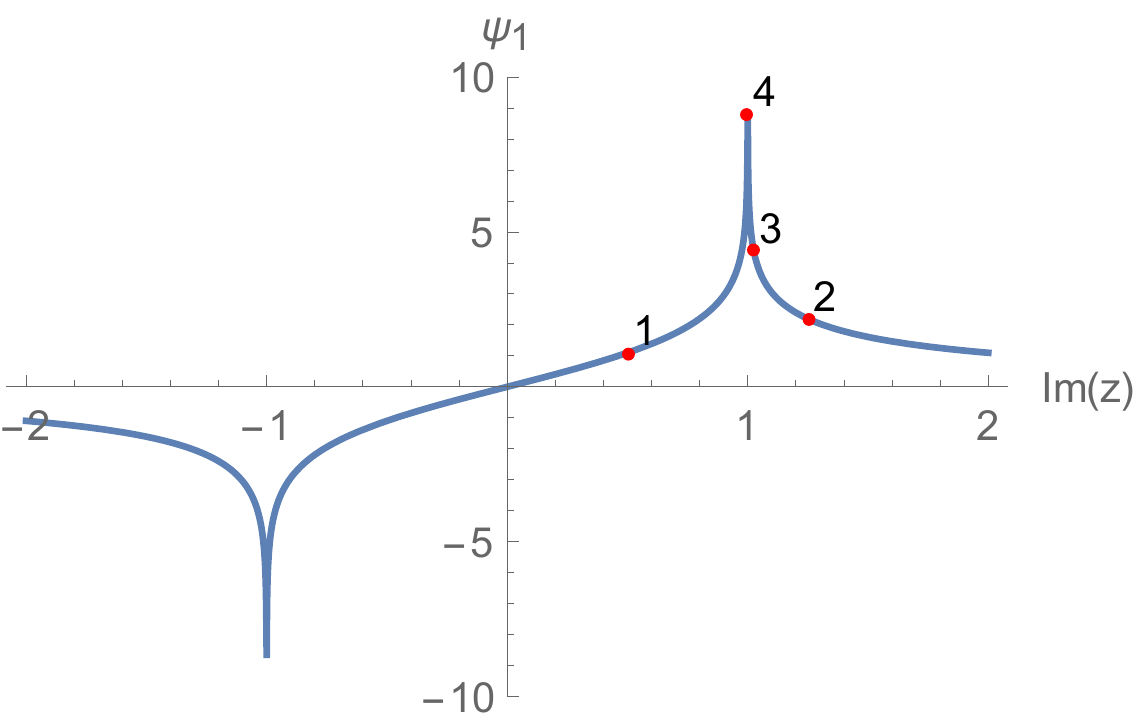}
\caption{\label{fig:newton_on_cubic_eigenfunctions}The eigenfunction $\psi_1$ of $\K_{N_c}$ for $a=1$, $b=0$, $c=1$, and $d=1$, with corresponding polynomial $f(z)=z^2+1$. Left: Plot on the complex plane, with contour lines labeled by the function values. Right: Eigenfunction evaluated on $[-2i,2i]$. The numbered, red dots show four evaluations of the Newton-Raphson method with starting point $z_1=0.5i$, $\psi_1(z_1)\approx 1.0986$. As expected, $\psi_1(z_{k+1})=2^{k}\psi_1(z_1)$ and the function diverges at the roots $\pm i$ of $f$.}
\end{figure}

In this section, we mostly use the conjugacy results from \cite{miller-2006} to construct eigenfunctions for the Koopman operator of the Newton-Raphson method for root finding.
We only consider polynomial functions of degree two in this analysis, and at the end give a hint of more complicated dynamics present for cubic polynomials.
The Newton-Raphson method to find roots of a function $f$ is defined as the iterative scheme
\begin{equation*}
x_{n+1}=x_n-[J{f}]^{-1}(x_n) {f(x_n)}=:N_f(x_n),
\end{equation*}
where $[Jf]_{ij}=\frac{\partial f_i}{\partial x_j}$ is the Jacobian matrix of $f$. In the example we consider here, the function $f:\mathbb{C}\to\mathbb{C}$ is a complex polynomial of degree two, for which the Newton-Raphson method simplifies to
\begin{equation*}
N_f(z)=z-\frac{f(z)}{f'(z)}.
\end{equation*}
Without loss of generality, we assume $f$ has the form $f(z)={a} z^2+b z + d$ with $a,b,d\in\mathbb{C}$, $a\neq 0$.
Any polynomial $f$ of degree two is conjugate to exactly one polynomial $g_c(z)=z^2+c$, $c\in\mathbb{C}$ (see \cite{miller-2006}), where the conjugacy is given through a linear transformation $h(z)=az+\frac{b}{2}$, such that
\begin{equation*}
f\circ h = h\circ g_c.
\end{equation*}
The constant $c$ is related to $a,b,d$ through $c=ad+\frac{b}{2}-\frac{b^2}{4}$.
The map $h$ scales and translates the argument so that $g_c$ has the desired form.
The conjugacy of $f$ and $g_c$ implies that we only need to analyze the Newton-Raphson method applied to functions of the form $g_c$, with a given $c$.
For this simpler form of $f$, the Newton-Raphson method becomes
\begin{equation*}
N_{c}(z)=z-\frac{z^2+c}{2z}=\frac{z^2-c}{2z}.
\end{equation*}
In the Koopman operator picture, the map $N_{c}$ defines the dynamical system in discrete time with an associated Koopman operator $\K_{N_c}$.
An important result in \cite{miller-2006} is that the map $N_c$ is conjugate to the polynomial $g_0(z)=z^2$ through the map 
\begin{equation*}
h_0(z)=\frac{z+i\sqrt{c}}{z-i\sqrt{c}},
\end{equation*}
which is a M\"obius transformation, and hence invertible on ${\mathbb{C}\cup\{\infty\}}$. Thus, eigenfunctions of $\K_{N_c}$ can be constructed through eigenfunctions of $\K_{g_0}$, which are defined through
\begin{equation*}
[\K_{g_0}\phi_k](z)=(\Kefunc_k\circ g_0)(z)=\Keval_k \Kefunc_k(z),\ k\in\mathbb{Z}.
\end{equation*}
Since $g_0(z)=z^2$, there are eigenfunctions $\phi_k$ of the form $\Kefunc_k(z)=\ln(|z|)^k$, associated to the eigenvalues $\Keval_k=2^{k}$.
Through the given conjugacy $h_0$, we have that 
\begin{equation}\label{eq:koopman eigenfunctions newton}
\psi_k(z)=(\phi_k\circ h_0)(z)=\ln\left(\left|\frac{z+i\sqrt{c}}{z-i\sqrt{c}}\right|\right)^k
\end{equation}
are eigenfunctions for $\K_{N_c}$, associated to the eigenvalues $\Keval_k$. Fig.~\ref{fig:newton_on_cubic_eigenfunctions} shows $\psi_1$ for $c=1$, i.e. for the polynomial $f(z)=z^2+1$. The roots at $\pm i$ are clearly visible as the points where the eigenfunction diverges.

\subsubsection{An example with chaotic behavior}
\begin{figure}[ht!]
\centering
\includegraphics[width=\figureHalfWidth]{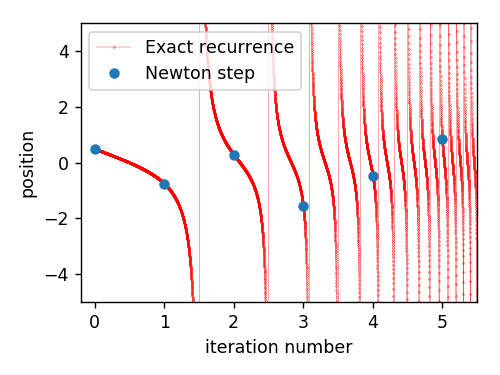}
\includegraphics[width=\figureHalfWidth]{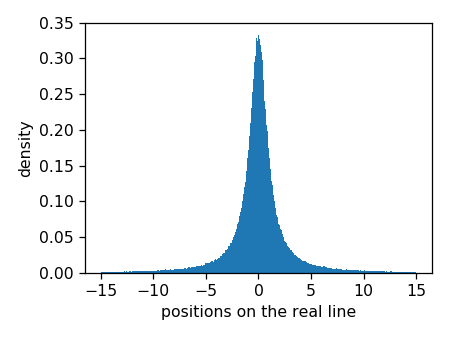}
\caption{\label{fig:newton_iterations_chaos}Left: The first five iterations of the Newton-Raphson method (blue dots) on the polynomial $f(z)=z^2+1$, with the initial condition $z(0)=1/2$. The exact recurrence relation is shown as red lines, evaluated over a continuous space, showing that the expansion and folding leads to faster and faster oscillations, a sign for chaotic behavior on the discrete iterations. Right: Invariant density over the real line, approximated from a single trajectory of length $n=10^6$, starting at $z_0=0.3$.}
\end{figure}

Consider the function $f(z)=z^2+1$ from the previous section, with its two roots at $\pm i$. Interpreting the Newton-Raphson method as a discrete dynamical system, the real line is the basin boundary separating the two basins of attraction on the complex plane. 
When started exactly on the real line, the Newton-Raphson method exhibits chaotic behavior, which we analyze in this section.
The Newton-Raphson method applied to $f$ is defined by the map
\begin{equation}\label{eq:newton on cubic}
N_f(z)=z-\frac{z^2+1}{2z}=\frac{z^2-1}{2z},
\end{equation}
so that $z_{n+1}=N_f(z_n)$ yields a new iterate, starting at $z_0\in\mathbb{C}$.
For $z_0\in\mathbb{C}/\mathbb{R}$, the limit points of the iterative scheme are the extremal points, where
\begin{equation*}\begin{matrix}
\lim_{n\to\infty}z_n&=&i&\Imag(z_0)>0,\\
\lim_{n\to\infty}z_n&=&-i&\Imag(z_0)<0.
\end{matrix}
\end{equation*}
For $z_0\in\mathbb{R}$, however, the scheme does not converge (see Fig.~\ref{fig:newton_iterations_chaos}) and gives rise to the chaotic recurrence (which can be solved explicitly in this case)
\begin{equation}\label{eq:newton on cubic true}
z(n)=-\cot\left(c_1 2^n\right),\ c_1=\arctan(-1/z_0).
\end{equation}
A similar discrete system with explicit solution is given by the logistic map at parameter $r=4$, with recurrence relation $y_{n+1}=r y_{n}(1-y_{n})$ and solution $y_n=\sin(c_0\, 2^n)^2$, $c_0=\arcsin\left(\sqrt{y_0}\right)$.
Accordingly, except for the constant function, the continuous eigenfunctions of the Koopman operator computed in the previous section are all zero (or infinity, for $k<0$) on the real line, see equation~(\ref{eq:koopman eigenfunctions newton}). For the map $N_f$, the real line is an invariant, repelling limit set of measure zero.

\begin{figure}[ht!]
\centering
\includegraphics[width=\figureFullWidth]{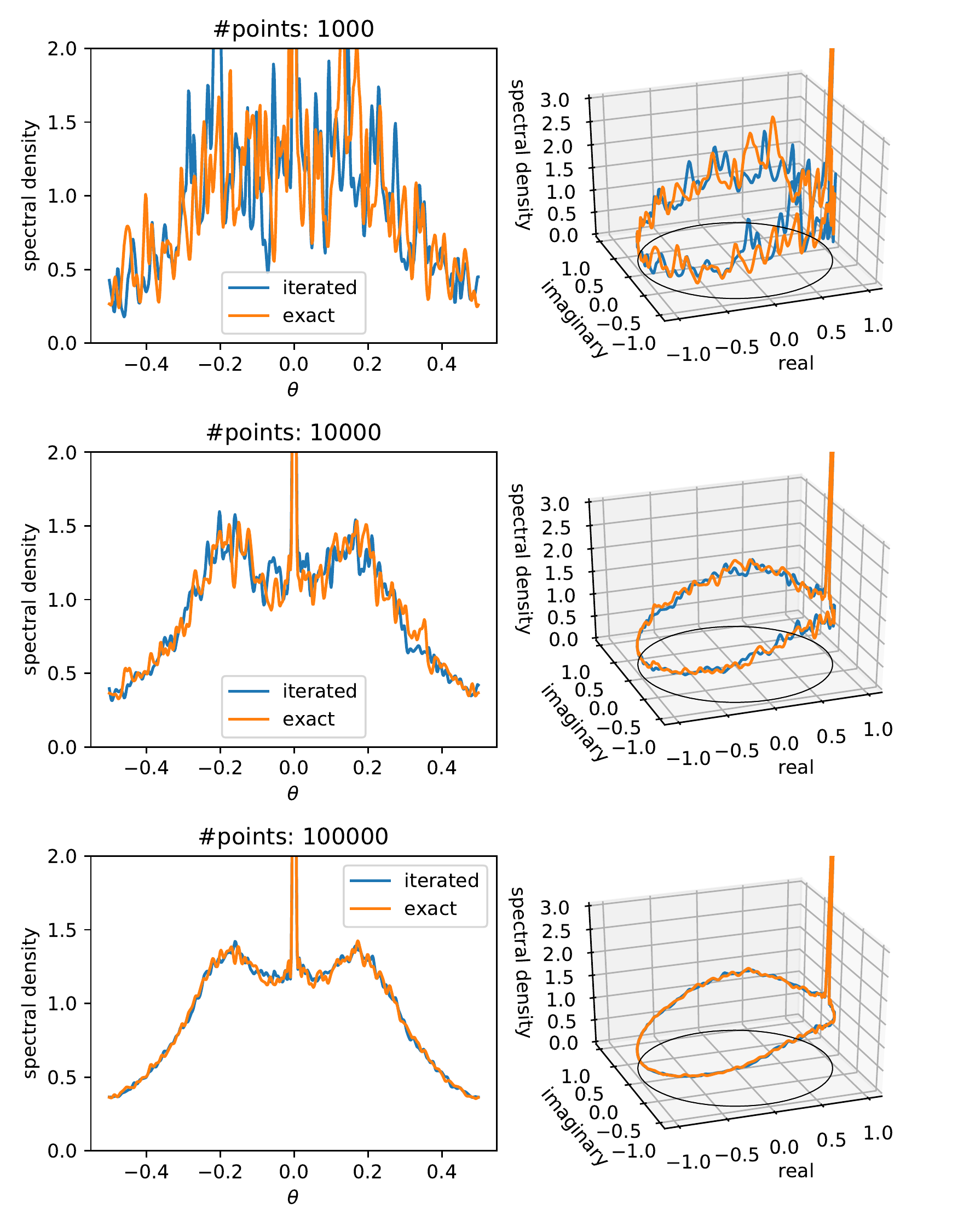}
\caption{\label{fig:newton_spectral_density}The spectral density on the unit circle $\mathbb{T}$, associated to the observable $g(z)=1+\exp(2\pi i z)$, is approximated with single trajectories of Newton's method with length $10^3$, $10^4$, and $10^5$ (rows) starting at $z_0=1/2+0i$. The point spectrum with eigenvalue at $1+0i$ associated to the constant function and the continuous spectrum on $\mathbb{T}/\{1\}$ can be distinguished. The blue curves are obtained from iterating the map $N_f$ with 16 digits of accuracy, the orange curves by using the exact formula and 35,000 digits of accuracy.}
\end{figure}

\review{The results so far have been constructed explicitly. Addressing the main challenge of the manuscript, we now want to demonstrate that numerical algorithms can also be used to analyze the behavior of the Newton-Raphson algorithm. To this end,
we employ} the numerical method by Korda, Putinar and Mezi\'{c}~\cite{korda-2018} (as described in \secref{approx continuous spectrum}) to approximate the continuous spectrum of the Koopman operator associated to the map $N_f$. \Figref{newton_spectral_density} shows how different numbers of points (i.e. lengths of trajectories) approximate the measure $\mu_g$ for $g(z)=1+\exp(2\pi i z)$, with $z\in\mathbb{R}$ from a trajectory of the Newton-Raphson method starting at $z_0=1/2$. The point spectrum is $\sigma_{pp}=\{1\}$ with eigenvalue $1$ associated to the constant function, as expected for a chaotic system. The absolutely continuous spectrum appears to be $\sigma_{c}=\mathbb{T}/\{1\}$, with no singularly continuous part.
The trajectories were computed with the iteration of $N_f$ (\eqnref{newton on cubic}) with 16 digits of accuracy, and with the explicit formula (\eqnref{newton on cubic true}) with 35,000 digits of accuracy, using the Wolfram Mathematica 11 software. The high accuracy is necessary to exactly represent the number $2^{100,000}\approx 10^{31,000}$ for the last point in the longest trajectory.

\subsubsection{Fractal eigenfunctions}
\begin{figure}[ht!]
\centering
\includegraphics[width=\figureFullWidth]{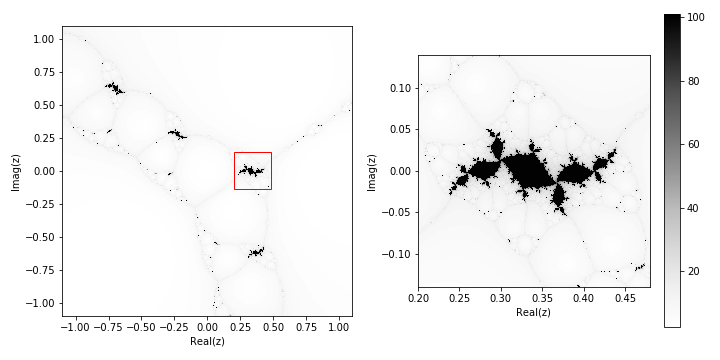}
\caption{\label{fig:newton_iterations_cubic}\review{Number of iterations $n$ (colorbar) until convergence up to $|N_f(z^*)-z^*|<0.01$ or $n\geq 100$}, for a cubic polynomial function $f(z)=(z+w)(z-w)(z-1)$ with $w=.589 + .605j$. The fractal structure of the number of iterates is visible in both panels. The right panel shows a zoomed in version of the left (marked with a box on the left).}
\end{figure}

\begin{figure}[ht!]
\centering
\includegraphics[width=0.45\textwidth]{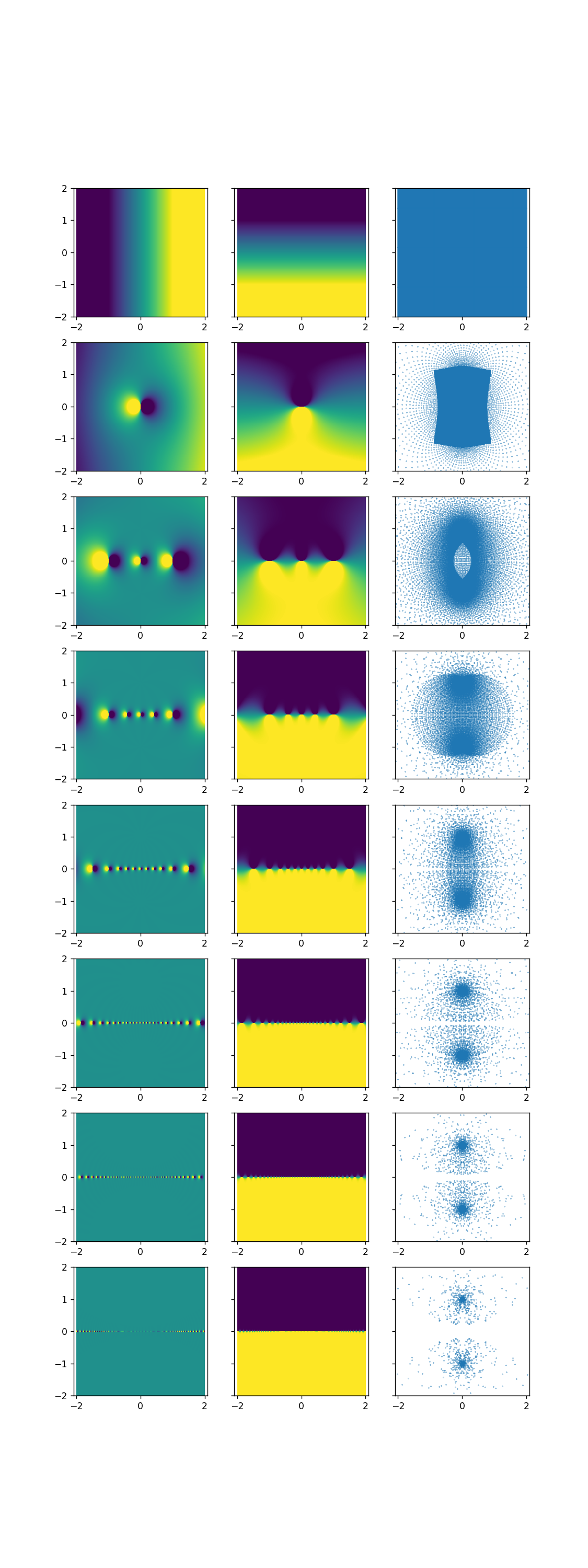}
\includegraphics[width=0.45\textwidth]{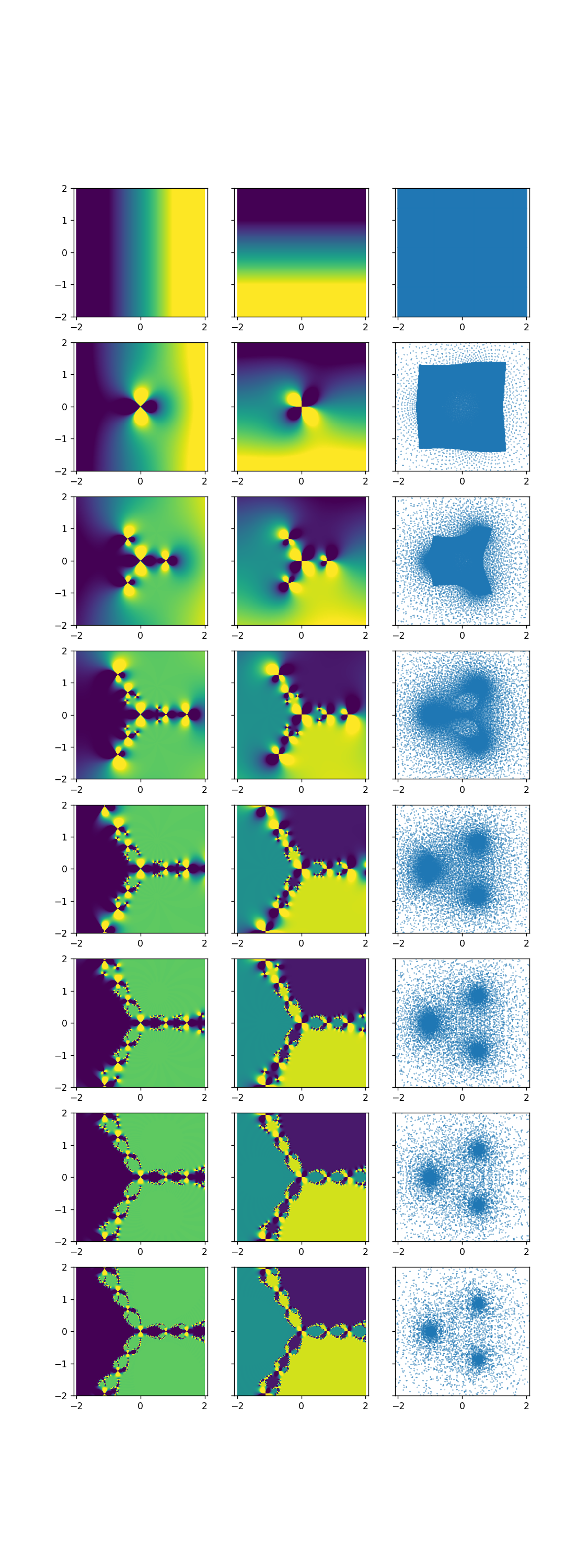}
\caption{\label{fig:newton_on_zsqrplusone}Iterations on the complex plane for a polynomial function of degree two (left) and degree three (right). The eigenfunctions at eigenvalue 1 are limits of the process. On the left, the eigenfunction has two separate values for the two basins of attraction. On the right, the eigenfunction has a persistent fractal structure.}
\end{figure}

For cubic polynomial functions, application of the {New\-ton}-{Raphson} method is known to exhibit an even richer structure compared to the quadratic polynomials described in the previous section.
We only give numerical results, inspired by \cite{miller-2006}.
\Figref{newton_iterations_cubic} shows the number of iterations (in color) of the Newton-Raphson method until the derivative of the objective function is smaller than $0.01$. 
\Figref{newton_on_zsqrplusone} shows eight iterations of the Newton\--Raphson method on a subset of points sampled on the complex plane, for a polynomial of degree two (left three columns), and degree three (right three columns). The convergence to the two (and three) roots is visible in the third and sixth column. The first/second and fourth/fifth columns show the real/imaginary parts of the final positions as functions on the initial positions, which makes the fractal structure for degree three (and higher) apparent.


\section{Conclusions}\label{sec:conclusions}
\review{Many modern algorithms (e.g. those associated with training neural networks) can difficult to analyze in a classical sense: they are complex, interconnected, can be deterministic or stochastic, act on high-dimensional states, and are often only accessible as input-output blax box systems. To address this challenge, we propose to employ the unifying view of the Koopman operator framework, with several numerical approximations readily available.
In a series of examples, we demonstrated that the Koopman operator provides such a unifying view of many important issues in the analysis of algorithms}: (spectral) convergence, algorithm acceleration, state space decomposition, high-dimensional state spaces, partial information, and generating processes for discrete algorithms. In the last section, we gave an outlook to chaotic behavior of algorithms and how the behavior of statistics of observables can be analyzed through the spectrum of the operator.
We discussed the possibility of using this approach to accelerate algorithms in high-dimensional embedding spaces whose long-term dynamics lie on low-dimensional manifolds. Another possibility is to use this approach to create surrogates of complicated ``black box" algorithms that are difficult to analyze mathematically; the data-driven surrogates obtained by sampling iterated algorithm states may provide useful insights in their nature and behavior. The Koopman operator has been shown to provide a convenient framework for constructing data-driven homeomorphisms between dynamical systems; it is an intriguing possibility that we can use this framework to realize homeomorphisms between different algorithms for solving the same problem.

\bibliographystyle{siamplain}
\bibliography{lit}

\appendix

\section{Computational experiments with partial information}\label{sec:partial information}
In most applications, the data will only partially cover the state space of an algorithm. In this section, we explore the behavior of the approximated operator spectrum on such partial domains.
The theoretical underpinning of this experiment is provided by the concept of ``open eigenfunctions''~\cite{mezic-2017}, i.e. eigenfunctions defined only on parts of the state space.
Not all the results from the computational experiments can be readily explained by the current theory, however, as we will detail below.

We study discrete gradient descent on the function $$f(x_1,x_2)=x_1^4-x_1^2+x_1/4+x_2^2.$$
\Figref{partial domain} illustrates the function values on $[-1,1]^2$ through contour lines (color), and the sample points used for the approximation of the operator (red area).
When using gradient descent to minimize $f$, we obtain a (discrete) dynamical system with two attracting steady states, and one saddle point in between them (see \figref{partial domain}).
The spectrum obtained by EDMD for the given sample domain (red) is shown on the right: all eigenvalues are inside the unit disk, indicating that the system behavior is purely attracting. This result is reasonable, because in the limit of infinite applications of the Koopman operator, all observable functions $\obs\in\mathcal{F}$ of the data set can be expressed by two piecewise constant eigenfunctions $\Kefunc_1^{(1)},\Kefunc^{(2)}_1$ associated to the eigenvalue $\Keval=1$:
\begin{eqnarray}\label{eq:limit behavior on partial domain}
\lim_{n\to\infty}\K^n \obs &=& a_1 \Kefunc^{(1)}_1 + a_2 \Kefunc^{(2)}_1;\ a_1,a_2\in\mathbb{C},\ g\in\mathcal{F},\\
\Kefunc^{(1)}_1(x)=\left\lbrace \begin{matrix}c_1^{(1)}&\text{if }x_1<x^{(s)}_1\\c_2^{(1)}&\text{if }x_1\geq x^{(s)}_1\end{matrix}\right.,&&\Kefunc^{(2)}_1(x)=\left\lbrace \begin{matrix}c_1^{(2)}&\text{if }x_1<x^{(s)}_1\\c_2^{(2)}&\text{if }x_1\geq x^{(s)}_1\end{matrix}\right.,
\end{eqnarray}
where $x^{(s)}_1$ is the $x_1$ coordinate of the saddle point between the two steady states.
The constants $a_1,a_2$ are determined by the value of $\obs$ on the two attracting steady states and the four constants $c_1^{(1)},c_1^{(2)},c_2^{(1)},c_2^{(2)}\in\mathbb{C}$ associated to the two eigenfunctions.
\Eqnref{limit behavior on partial domain} illustrates that, in the function space $\mathcal{F}$, there is a four-dimensional (real and imaginary parts of two complex numbers $a_1,a_2$), attracting subspace that is the limit set of the system $\obs_{n+1}=\K \obs_n$.

If we shift the domain to the right (\figref{partial domain max eigenvalue}, C), the left attracting steady state leaves the sampling domain. This causes trajectories starting in the sampling domain to leave it after some time, which typically indicates a repelling set---explaining the numerically obtained eigenvalues outside the unit disk in parts A and C of \figref{partial domain max eigenvalue}.
Note that since many trajectories leave the domain after a finite number $N$ of iterations (depending on the initial state $x_n$), the dynamical system is not defined for iteration numbers $n>N$. For points on the boundary of the domain, $N=0$, so the flow map at these points is only defined for iterations backwards in time.
Since there is a saddle inside the domain, the system is also not defined for some states $x\in X$ on the boundary that would leave the domain backwards in time (iteration numbers $n<0$).
This poses problems with the definition of the Koopman operator family $\{\K^n\}$, since it acts on functions defined on all of $X$---and for all $n\neq 0$, some elements leave $X$, so that the flow map on the entire set $X$ is not defined for $n\neq 0$.
If all vectors on the boundary of the data domain pointed inward (outward), the flow would exist for all forward (backward) time, and the problem would be related to inflowing (overflowing) invariant manifolds~\cite{eldering-2018}.
Eigenvalues with a real part larger than one indicate unstable behavior, indicating the existence of unstable nodes or saddles in the data set.
\begin{figure}[htp]
\centering
\includegraphics[width=\figureFullWidth]{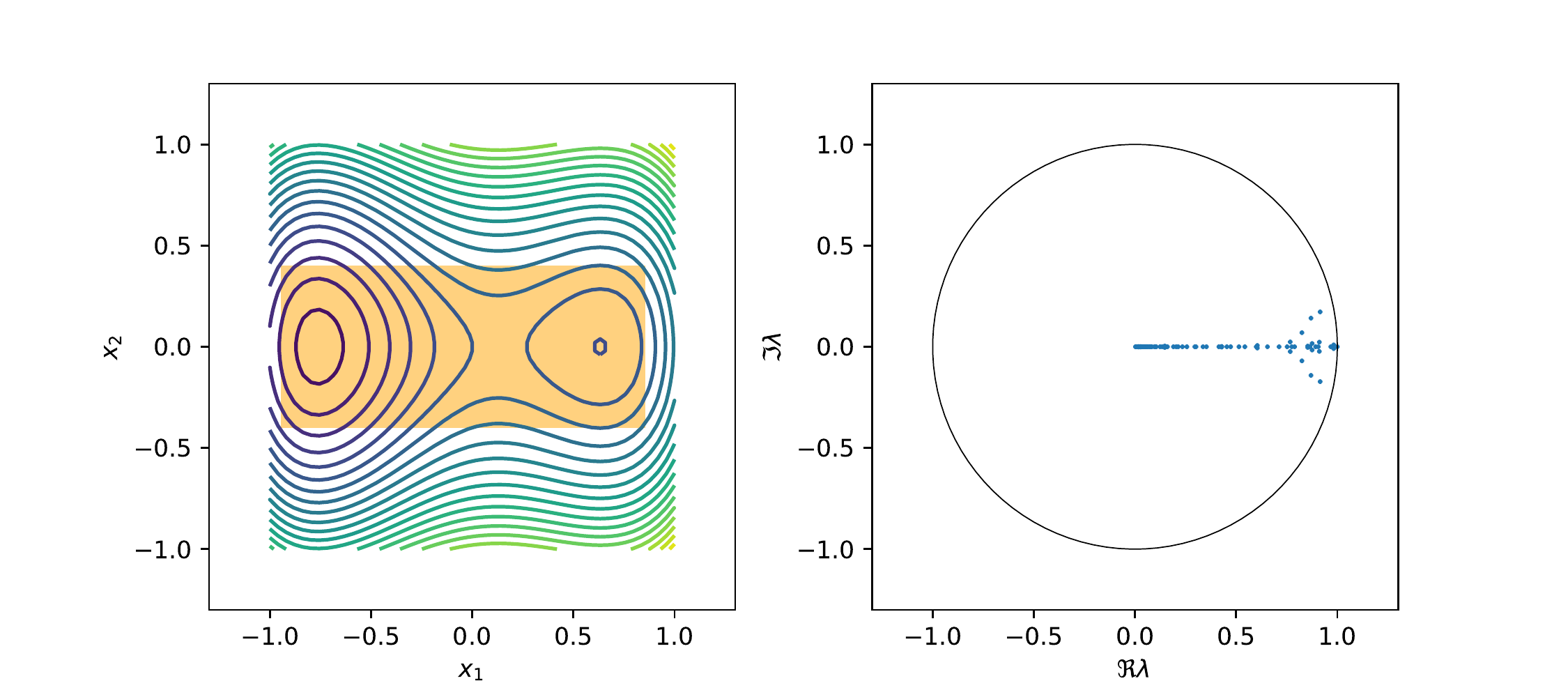}
\caption{\label{fig:partial domain}Partial sampling (orange rectangle) of a two-dimensional domain. The objective function is shown as a countour plot. The two attracting steady states are both inside the rectangle, and the spectrum is fully inside the unit disk (right plot).}
\end{figure}
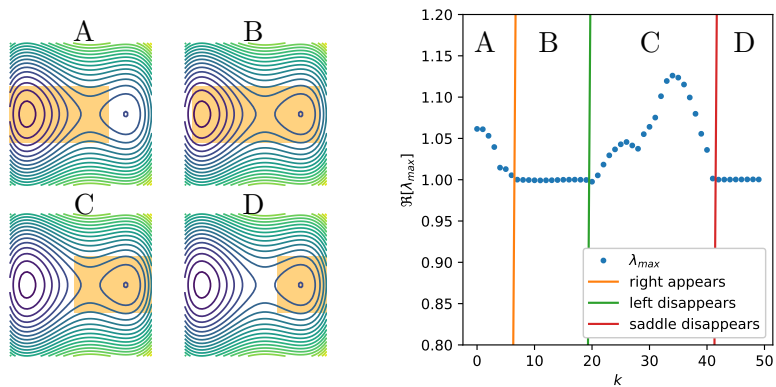
\begin{figure}[htp]
\centering
\begin{tikzpicture}
\node[inner sep=0pt] (figure) at (0,0)
   {\includegraphics[width=5.5cm]{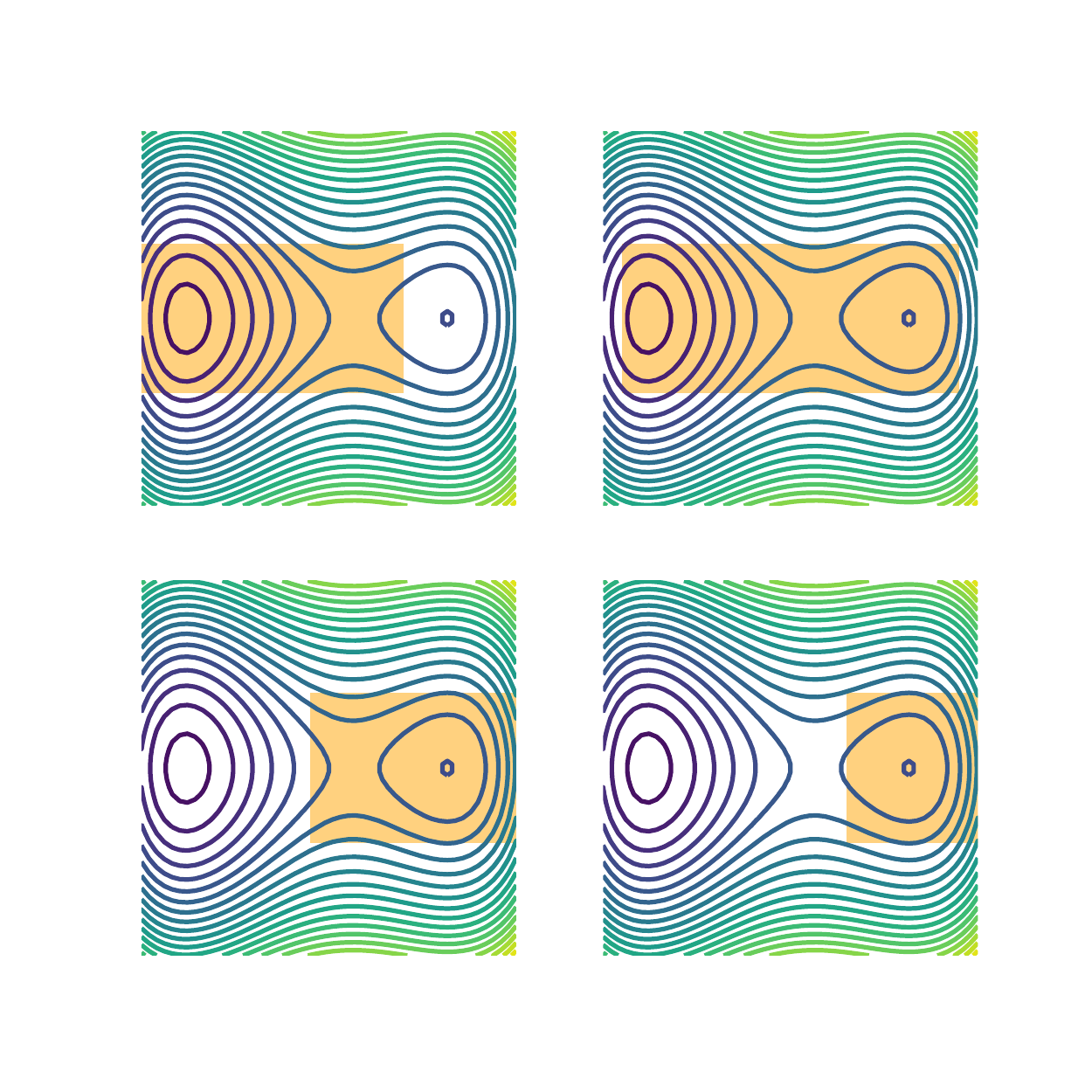}};
\node[inner sep=0pt] (a) at (-1.05,2.25){A};
\node[inner sep=0pt] (b) at (1.2,2.25){B};
\node[inner sep=0pt] (c) at (-1.05,-0.05){C};
\node[inner sep=0pt] (c) at (1.2,-0.05){D};
\end{tikzpicture}
\begin{tikzpicture}
\node[inner sep=0pt] (figure) at (2,0)
   {\includegraphics[width=5.5cm]{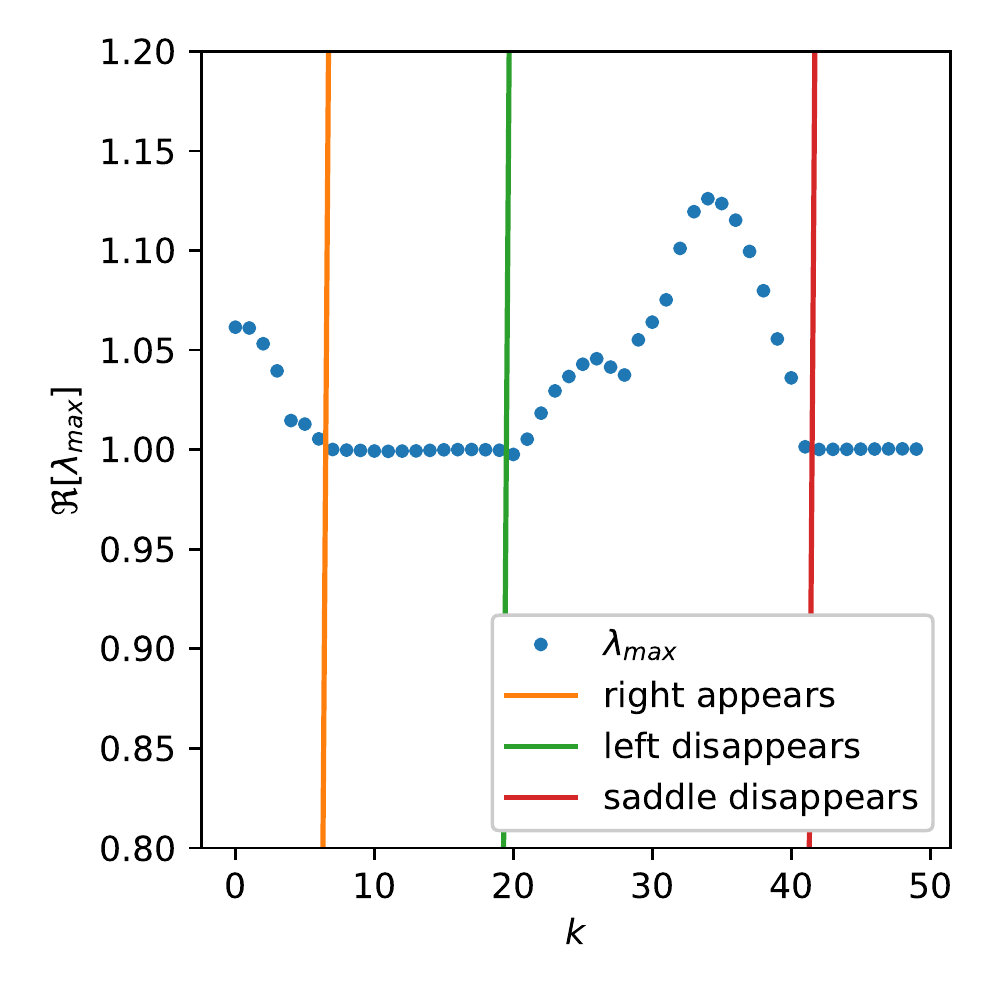}};
\node[inner sep=0pt] (a) at (.65,2.1){A};
\node[inner sep=0pt] (b) at (1.5,2.1){B};
\node[inner sep=0pt] (c) at (2.85,2.1){C};
\node[inner sep=0pt] (c) at (4.1,2.1){D};
\end{tikzpicture}
\caption{\label{fig:partial domain max eigenvalue}At the beginning, the right steady state is not inside the data domain (orange rectangle, A). When both steady states are contained in the data, the maximum of the real part of all eigenvalues is close to one (B). As the available data set moves to the right, the left steady state is no longer available to the approximation (C). When the rectangle is moved even further (D), the saddle between the two steady states also disappears and only one, attracting steady state is contained in the data.}
\end{figure}

\end{document}